\documentclass[11pt]{amsart}
\usepackage{graphicx, amscd, amsfonts, amsthm, amssymb}
\usepackage[all]{xy}
\title{Homotopy types of Diffeomorphism groups \\
of noncompact 2-manifolds}
\author[T.Yagasaki]{Tatsuhiko Yagasaki}
\subjclass[2000]{57S05, 58D05, 57N20, 58A10}
\keywords{Transformation groups, Topological groups, Diffeomorphism groups, 
Volume-preserving diffeomorphisms, End charge homomorphism, 
Infinite-dimensional manifolds, $\sigma$-compact manifold, Surfaces}
\address{Division of Mathematics, 
Graduate School of Science and Technology, 
Kyoto Institute of Technology,  
Matsugasaki, Sakyoku, Kyoto 606-8585, Japan}
\email{yagasaki@kit.ac.jp}

\newtheorem{theorem}{Theorem}[section]
\newtheorem{proposition}{Proposition}[section] 
\newtheorem{corollary}{Corollary}[section] 
\newtheorem{lemma}{Lemma}[section]

\theoremstyle{definition}
\newtheorem{defn}{Definition}[section]
\newtheorem{remark}{Remark}[section]

\newtheorem*{assumption}{Assumption (A)} 
\newtheorem*{condition}{Condition (C)}

\def \cal {\mathcal}
\def \phi {\varphi}

\begin{document}
\baselineskip 6 mm

\thispagestyle{empty}
\maketitle
\begin{abstract}
Suppose $M$ is a noncompact connected smooth 2-manifold without boundary and 
let ${\cal D}(M)_0$ denote the identity component of the diffeomorphism group of $M$ with the compact-open $C^\infty$-topology. 
In this paper we investigate the topological type of ${\cal D}(M)_0$ and 
show that ${\cal D}(M)_0$ is a topological $\ell_2$-manifold and it has the homotopy type of the circle if $M$ is the plane, 
the open annulus or the open M\"obius band, and 
it is contractible in all other cases. 
When $M$ admits a volume form $\omega$, 
we also discuss the topological type of the group of 
$\omega$-preserving diffeomorphisms of $M$. 
To obtain these results we study some 
fundamental 
properties of transformation groups on noncompact spaces endowed with weak topology. 
\end{abstract}

%%%%%%%%%%%%%%%%%%%%%%%%%%%%

\section{Introduction}

The purpose of this paper is the investigation of topological properties of the diffeomorphism groups 
of noncompact smooth 2-manifolds endowed with the {\it compact-open} $C^\infty$-topology. 
When $M$ is a closed smooth $n$-manifold, 
the diffeomorphism group ${\cal D}(M)$ with the compact-open $C^\infty$-topology is a smooth Fr\'echet manifold \cite[Section I.4]{Ham}, and 
for $n = 2$, S.~Smale \cite{Sm} and C.~J.~Earle and J.~Eell \cite{EE} classified the homotopy type of the identity component ${\cal D}(M)_0$. 

In the $C^0$-category, for any compact 2-manifold $M$, 
the homeomorphism group ${\cal H}(M)$ with the compact-open topology is a topological $\ell_2$-manifold \cite{DT, Ge, LM, To}, and 
M.~E.~Hamstrom \cite{Ha} classified the homotopy type of the identity component ${\cal H}(M)_0$ (cf. \cite{Sc} for PL-case). 
In \cite{Ya2} we have shown that ${\cal H}(M)_0$ is an $\ell_2$-manifold even if $M$ is a noncompact connected 2-manifold.
We also classified its homotopy type and showed that ${\cal H}(M)_0$ is contractible except a few cases. 

In \cite{BMSY} we studied topological types of 
transformation groups on noncompact spaces endowed with strong topology. 
This formulation was intended for an application to 
homeomorphism groups and diffeomorphism groups of noncompact manifolds endowed with the Whitney topology. 

In this article we formulate the notion of weak topology for transformation groups on noncompact spaces. This notion corresponds with the compact-open topology. 
We see that the main arguments in \cite{Ya1, Ya2, Ya3} well extend 
to transformation groups with weak topology (cf. Theorem~\ref{thm_criterion}) and  
these results can be well applied to the diffeomorphism groups of noncompact 2-manifolds. 

Suppose $M$ is a smooth $n$-manifold and $X$ is a closed subset of $M$. 
For $r = 1, 2, \cdots, \infty$ we denote by ${\cal D}^{\, r}_X(M)$ the group of $C^r$-diffeomorphisms $h$ of $M$ onto itself with $h|_X = id_X$, 
endowed with the compact-open $C^r$-topology 
\cite[CH.2 Section 1]{Hi}, 
and by ${\cal D}_X^{\, r}(M)_0$ the identity connected component of ${\cal D}^{\, r}_X(M)$. 
By a compact smooth submanifold of $M$ we mean 
the union of a disjoint family 
of a closed smooth $k$-submanifold of $M$ for $k = 0, 1, \cdots, n-1$ and 
a compact smooth $n$-submanifold of $M$.  
The following is the main result of this paper. 

\begin{theorem}\label{thm_l2}
Suppose $M$ is a noncompact connected smooth $2$-manifold without boundary and $X$ is a compact smooth submanifold of $M$. Then the following hold. 
\begin{itemize}
\item[{\rm (1)}] ${\cal D}^{\, r}_X(M)_0$ is a topological $\ell_2$-manifold. 
\item[{\rm (2)}] 
\begin{itemize}
\item[{\rm (i)}\,] ${\cal D}^{\, r}_X(M)_0 \simeq {\Bbb S}^1$ if $(M, X) =$\hspace{-2mm} 
\begin{tabular}[t]{l}
$($a plane, $\emptyset)$, $($a plane, 1\,pt$)$, $($an open M\"obius band, $\emptyset)$ \\[1mm]
\ \ or $($an open annulus, $\emptyset)$. 
\end{tabular} 
\item[{\rm (ii)}] ${\cal D}^{\, r}_X(M)_0 \simeq \ast$ in all other cases.
\end{itemize} 
\end{itemize}
\end{theorem}

Note that any separable infinite-dimensional Fr\'echet space 
is homeomorphic to the separable Hilbert space $\ell_2 \equiv \{ (x_n) \in {\Bbb R}^\infty : \sum_n x_n^2 < \infty \}$ \cite[Chapter VI, Theorem 5.2]{BP}. 
A topological $\ell_2$-manifold is a separable metrizable space which is locally homeomorphic to $\ell_2$. 
Since topological types of $\ell_2$-manifolds are classified by their homotopy types, 
Theorem 1.1\,(2) implies that ${\cal D}^{\, r}_X(M)_0 \cong {\Bbb S}^1 \times \ell_2$ in the case (i) 
and ${\cal D}^{\, r}_X(M)_0 \cong \ell_2$ in the case (ii). 

Let ${\cal H}_X(M)_0$ denote the identity connected component of the group of homeomorphisms $h$ of $M$ onto itself with $h|_X = id_X$,  
endowed with the compact-open $C^0$-topology. 
The comparison of the homotopy types of ${\cal D}^{\, r}_X(M)_0$ and ${\cal H}_X(M)_0$ (\cite{Ya2}) implies the following conclusion:

\begin{corollary}
Suppose $X$ is a compact smooth submanifold of $M$. 
Then the inclusion ${\cal D}^{\, r}_X(M)_0 \subset {\cal H}_X(M)_0$ is a homotopy equivalence. 
\end{corollary}

For the subgroup of diffeomorphisms with compact supports, we have the following consequences. 
Let ${\cal D}_X^{\, r}(M)^c$ denote the subgroup of ${\cal D}^{\, r}_X(M)$ consisting of diffeomorphisms with compact supports, 
and let ${\cal D}_X^{\, r}(M)^c_0$ denote the identity connected component of ${\cal D}_X^{\, r}(M)^c$. 
We can also consider the subgroup of diffeomorphisms which are isotopic to $id_M$ by isotopies with compact supports.
Let ${\cal D}_X^{\, r}(M)^{c\,\ast}_0$ denote the subgroup of ${\cal D}_X^{\, r}(M)^c_0$ consisting of 
$h \in {\cal D}_X^{\, r}(M)^c$ which admits an ambient $C^r$-isotopy $h_t : M \to M$ rel $X$ such that $h_0 = h$, $h_1 = id_M$ and 
$h_t$ ($0 \leq t \leq 1$) have supports in a common compact subset of $M$.  

We say that a subspace $A$ of a space $X$ is homotopy dense (or has the homotopy negligible complement) in $X$ 
if there exists a homotopy $\phi_t : X \to X$ such that $\phi_0 = id_X$ and $\phi_t(X) \subset A$ ($0 < t \leq 1$). 
In this case, the inclusion $A \subset X$ is a (controlled) homotopy equivalence, and 
when $X$ is metrizable, $X$ is an ANR iff $A$ is an ANR (cf.~\S2.4). 

\begin{theorem}\label{thm_HD}
Suppose $M$ is a noncompact connected smooth $2$-manifold without boundary and $X$ is a compact smooth submanifold of $M$. Then ${{\cal D}_X^{\, r}(M)^c_0}^\ast$ is homotopy dense in ${\cal D}^{\, r}_X(M)_0$ 
\end{theorem} 

\begin{corollary}\label{cor_HD}
{\rm (1)} ${\cal D}_X^{\, r}(M)^c_0$ and ${{\cal D}_X^{\, r}(M)^c_0}^\ast$ are ANR's.  
\begin{itemize}
\item[{\rm (2)}] The inclusions \ ${{\cal D}_X^{\, r}(M)^c_0}^\ast \subset {\cal D}_X^{\, r}(M)^c_0 \subset {\cal D}^{\, r}_X(M)_0$ \ 
are homotopy equivalences.
\end{itemize}
\end{corollary}

We notice that the subgroup ${{\cal D}_X^{\, r}(M)^c_0}^\ast$ coincides with ${\cal D}_X^{\, r}(M)^c_0$ except some specific cases. 

\begin{proposition}\label{prop_c*=c} Suppose $M$ is a noncompact connected smooth $2$-manifold without boundary and 
$X$ is a compact smooth 2-submanifold of $M$. Then 
${{\cal D}_X^{\, r}(M)^c_0}^\ast = {\cal D}_X^{\, r}(M)^c_0$ iff 
{\rm (a)} $M$ has no product end or {\rm (b)} $(M, X) = (\text{a plane}, \emptyset)$ or $($an open M\"obius band, $\emptyset)$. 
\end{proposition}

For $n$-manifolds of finite type we can deduce the following conclusion. 

\begin{proposition}\label{prop_finite-type}
If $M = {\rm Int}\,N$ for some compact smooth $n$-manifold $N$ with nonempty boundary and $X$ is a compact smooth submanifold of $M$, then 
${\cal D}^{\, r}_X(M)_0$ is an $\ell_2$-manifold. 
\end{proposition}

This proposition is not so obvious, because 
for the subgroup 
\[ \mbox{${\cal D}^{\, r}_{\partial N}(N)' = \{ h \in {\cal D}^{\, r}_{\partial N}(N) \mid f = id$ on a neighborhood of $\partial N \}$} \]  
the restriction map \ ${\cal D}^{\, r}_{\partial N}(N)' \to {\cal D}^r({\rm Int}\,N)^c$ \ is a continuous bijection, but not a homeomorphism. 

At this point it is important to compare the above results on the compact-open $C^\infty$-topology  %in this article 
with those on the Whitney $C^\infty$-topology in \cite{BMSY}.  
Suppose $M$ is a noncompact connected smooth $n$-manifold without boundary. 
Let ${\cal D}^\infty(M)_w$ denote the group ${\cal D}^\infty(M)$ endowed with the Whitney $C^\infty$-topology.  In \cite{BMSY} we have shown that 
${\cal D}^\infty(M)_w$ is locally homeomorphic $\square^\omega l_2$ (the countable box product of $\ell_2$), while 
${\cal D}^\infty(M)_w^c$ is an ${\Bbb R}^\infty\times l_2$-manifold. 
Here, ${\Bbb R}^\infty$ is the direct limit of the tower ${\Bbb R}^1\subset {\Bbb R}^2\subset {\Bbb R}^3\subset\cdots.$
Since the Whitney $C^\infty$-topology is so strong, 
it is seen that $({\cal D}^\infty(M)_w)_0 = ({\cal D}^\infty(M)^c_w)_0$ and 
any compact subset of ${\cal D}^\infty(M)^c_w$ has a common compact support. 
When $n = 2$, 
the difference between ${\cal D}^\infty(M)_0$ and $({\cal D}^\infty(M)_w)_0$ is summarized as follows. 

\begin{remark} If $M$ is a noncompact connected smooth $2$-manifold without boundary, then 
\vspace{1mm} 
\begin{itemize}
\item[(1)] ${\cal D}^\infty(M)_0 \cong \left\{ \hspace{-1mm} 
\begin{array}[c]{cl}
{\Bbb S}^1 \times \ell_2 & \text{if $M=$ a plane, an open M\"obius band, an open annulus,} \\[1mm]
\ell_2 & \text{in all other cases,} 
\end{array}\right.$ 
\vskip 2mm 
\item[] the inclusion ${\cal D}^\infty(M)^{c\,\ast}_0 \subset {\cal D}^\infty(M)_0$ is a homotopy equivalence, 
\vskip 1mm   
\item[(2)]  $({\cal D}^\infty(M)_w)_0 = ({\cal D}^\infty(M)^c_w)_0 \cong {\Bbb R}^\infty\times l_2$ \quad ($({\cal D}^\infty(M)_w)_0 = {\cal D}^\infty(M)^{c\,\ast}_0$ as sets).  
\end{itemize}
\end{remark} 

Note that any loop in $({\cal D}^\infty(M)_w)_0$ has a common compact support so that it is inessential, while 
%in the case where ${\cal D}^\infty(M)^{c\,\ast}_0 \simeq {\Bbb S}^1$ 
a loop $(h_t)_t$ in ${\cal D}^\infty(M)^{c\,\ast}_0$ may not have a common compact support though each $h_t$ admits an isotopy to $id_M$ with compact support. 
In \cite[Example 3.1\,(2)]{Ya3}) a homotopy equivalence 
$f : {\Bbb S}^1 \simeq {\cal H}({\Bbb R}^2)^{c\,\ast}_0$ is constructed explicitly. 
This fact can be regarded as a sort of pathology of the compact-open topology. 

Next we discuss the groups of volume-preserving diffeomorphisms of noncompact 2-manifolds. 
Suppose $M$ is a connected oriented smooth $n$-manifold possibly with boundary and 
$\omega$ is a positive volume form on $M$. 
Let ${\cal D}^\infty(M; \omega)$ denote the subgroup of ${\cal D}^\infty(M)$ 
consisting of $\omega$-preserving diffeomorphisms of $M$ (endowed with the compact-open $C^\infty$-topology)
and let ${\cal D}^\infty(M; \omega)_0$ denote the identity connected component of ${\cal D}^\infty(M; \omega)$. 
In \cite[Corollary 1.1]{Ya4} we have shown that the group ${\cal D}^\infty(M)_0$ has a factorization 
$$({\cal D}^\infty(M)_0, {\cal D}^\infty(M; \omega)_0) \cong 
({\cal V}^+(M; \omega(M), E_M^\omega)_{ew}, \{ \omega \}) \times {\cal D}^\infty(M; \omega)_0,$$
where ${\cal V}^+(M; \omega(M), E_M^\omega)$ is a convex space of positive volume forms on $M$ endowed with the finite-ends weak $C^\infty$-topology $ew$ (cf. Section 6).
This means that the subgroup ${\cal D}^\infty(M; \omega)_0$ is a strong deformation retract (SDR) of ${\cal D}^\infty(M)_0$, 
Hence Theorem~\ref{thm_l2} yields the following consequences. 

\begin{theorem}\label{thm_vol-pre} Suppose $M$ is a noncompact connected orientable smooth 2-manifold without boundary and $\omega$ is a volume form on $M$. 
Then the following hold. 
\begin{itemize}
\item[{\rm (1)}] ${\cal D}^\infty(M; \omega)_0$ is a topological $\ell_2$-manifold and it is 
a SDR of ${\cal D}^\infty(M)_0$. 
\item[{\rm (2)}] 
\begin{itemize}
\item[{\rm (i)}\,] ${\cal D}^\infty(M; \omega)_0 \simeq {\Bbb S}^1$ if $M =$ a plane or an open annulus. 
\item[{\rm (ii)}] ${\cal D}^\infty(M; \omega)_0 \simeq \ast$ in all other cases.
\end{itemize} 
\end{itemize}
\end{theorem}

S.~R.~Alpern and V.\,S.\,Prasad \cite{AP} introduced 
the end charge homomorphism 
$$c^\omega_0 : {\cal D}^\infty(M; \omega)_0 \to {\cal S}(M; \omega).$$ 
Here ${\cal S}(M; \omega)$ is the topological linear space of end charges of $(M, \omega)$ and 
for each $h \in {\cal D}^\infty(M; \omega)_0$
the end charge $c^\omega_0(h)$ measures volume transfer toward ends of $M$ under $h$ (cf. Section 6). 
In \cite[Corollary 1.2]{Ya4} we have shown that the homomorphism $c^\omega$ has a continuous (non-homomorphic) section.  
This induces the factorizations 
$$({\cal D}^\infty(M; \omega)_0, {\rm ker}\,c^\omega_0) 
\cong ({\cal S}(M; \omega), 0) \times {\rm ker}\,c^\omega_0$$
and  implies that the kernel  
${\rm ker}\,c^\omega_0$ is a SDR of 
${\cal D}^\infty(M; \omega)_0$. 

Let ${\cal D}^\infty(M; \omega)^c$  
denote the subgroup of ${\cal D}^\infty(M; \omega)$
consisting of $\omega$-preserving diffeomorphisms of $M$ with compact support and let ${\cal D}^\infty(M; \omega)^c_0$ denote 
the identity connected component of ${\cal D}^\infty(M; \omega)^c$. 
The subgroup ${\cal D}^\infty(M; \omega)^{c\,\ast}_0$ of ${\cal D}^\infty(M; \omega)^c_0$ is defined by 
$${\cal D}^\infty(M; \omega)^{c\,\ast}_0 = \big\{ h \in {\cal D}^\infty(M; \omega) \mid h \in {\cal D}^\infty_{M - K}(M; \omega)_0 \text{ for some compact subset $K$ of $M$}\big\}.$$
As another application of the results in Section 3, in $n = 2$ 
we can deduce the topological relations among the subgroups 
$${\cal D}^\infty(M; \omega)^{c\,\ast}_0 \ \subset \ {\cal D}^\infty(M; \omega)^c_0 \ \subset \ {\rm ker}\,c^\omega_0.$$

\begin{theorem}\label{thm_vol-pre_hd} 
Suppose $M$ is a noncompact connected orientable smooth 2-manifold without boundary and $\omega$ is a volume form on $M$. 
Then the following hold. 
\begin{itemize}
\item[{\rm (1)}] ${\rm ker}\,c^\omega_0$ is an $\ell_2$-manifold and it is 
a SDR of ${\cal D}^\infty(M; \omega)_0$. 
\item[{\rm (2)}] ${\cal D}^\infty(M; \omega)^{c\,\ast}_0$ is homotopy dense in ${\rm ker}\,c^\omega_0$. Therefore, 
\begin{itemize}
\item[(i)\,] ${\cal D}^\infty(M; \omega)^{c\,\ast}_0$ and ${\cal D}^\infty(M; \omega)^c_0$ are ANR's, and  
\item[(ii)] the inclusions \ ${\cal D}^\infty(M; \omega)^{c\,\ast}_0 \ \subset \ {\cal D}^\infty(M; \omega)^c_0 \ \subset \ {\cal D}^\infty(M; \omega)_0$ \ are homotopy equivalences. 
\end{itemize}
\end{itemize}
\end{theorem}

This paper is organized as follows. 
Section 2 contains some fundamental facts on ANR's and $\ell_2$-manifolds. 
Section 3 contains main arguments in this article. 
Here, we investigate some fundamental topological properties of transformation groups on noncompact spaces endowed with weak topology. 
In Section 4 these results are applied to the diffeomorphism groups of noncompact manifolds with the compact-open $C^r$-topology. 
Theorems 1.1 and 1.2 are proved in Section 5. 
In Section 6 we discuss the groups of volume-preserving diffeomorphisms of noncompact manifolds. 

%%%%%%%%%%%%%%%%%%%%%

\section{Basic properties of ANR's and $\ell_2$-manifolds} 

In Sections 3\,-\,6 we see that the ANR-property of diffeomorphism groups and embedding spaces is especially important 
to investigate the topology of diffeomorphism groups of noncompact manifolds.
In this section we recall basic properties of ANR's. We refer to \cite{Hu, Pa2} for the theory of ANR's. 
Throughout the paper we assume that spaces are separable and metrizable and maps are continuous (otherwise specified). 

A metrizable space $X$ is called an {\em ANR} (absolute neighborhood retract) for metrizable spaces 
if any map $f : B \to X$ from a closed subset $B$ of a metrizable space $Y$ admits an extension to a neighborhood $U$ of $B$ in $Y$. 
If we can always take $U = Y$, then $X$ is called an {\em AR}. An AR is exactly a contractible ANR.  
It is well known that $X$ is an AR (an ANR) iff it is a retract of (an open subset of) a normed space.  
Any ANR has a homotopy type of CW-complex. 

We will apply the following criterion of ANR's \cite{Hu}:

\begin{lemma}\label{lem_ANR} 
{\rm (1)} Any retract of of an AR $($an open subset of an ANR$)$ is an AR $($an ANR$)$. 
\begin{itemize}
\item[{\rm (2)}] A metrizable space $X$ is an ANR iff each point of $X$ has an ANR neighborhood in $X$. 
\item[{\rm (3)}] If $\displaystyle X = \cup_{i = 1}^\infty U_i$, $U_i$ is open in $X$, $U_i \subset U_{i+1}$ and each $U_i$ is an AR, then $X$ is an AR. 
\item[{\rm (4)}] $X \times Y$ is a nonempty ANR iff $X$ and $Y$ are nonempty ANR's. 
In a fiber bundle, the total space is an ANR iff both the base space and the fiber are ANR's. 
\item[{\rm (5)}] A metric space $X$ is an ANR iff 
for any $\varepsilon > 0$ there is an ANR $Y$ and maps $f : X \to Y$ and $g : Y \to X$ such that $gf$ is $\varepsilon$-homotopic to $id_X$ \cite{Han}.  
\end{itemize}
\end{lemma} 

In the statement (3) the space $X$ is an ANR and has the trivial homotopy groups, thus it is contractible.  
Since any Fr\'echet space is an AR, by (2) any Fr\'echet manifold (an $\ell_2$-manifold) is an ANR. 
 
In a fiber bundle, if the base space is contractible and paracompact, then 
this bundle is trivial. A principal bundle is trivial iff it admits a section. 
We use the following fact on principal bundles with AR fibers. 

\begin{lemma}\label{lem_bdle_AR} 
If a principal bundle $p : E \to B$ has an AR fiber and a metrizable base space $B$, then, 
\begin{itemize} 
\item[{\rm (i)}\,] the bundle $p$ admits a section $s$ and so it is a trivial bundle  and 
\item[{\rm (ii)}] the map $sp : E \to E$ is $p$-fiber-preserving homotopic to $id_E$ and so the map $p$ is a homotopy equivalence with a homotopy inverse $s$. 
\end{itemize} 
\end{lemma} 

The notion of homotopy denseness (or homotopy negligibility) has been defined in \S1. 
A subspace $A$ of a space $X$ is homotopy dense (HD) in $X$ 
if there exists a homotopy $\phi_t : X \to X$ such that $\phi_0 = id_X$ and $\phi_t(X) \subset A$ ($0 < t \leq 1$) (\cite{To}). 
The homotopy $\phi_t$ is called an absorbing homotopy.
The map $\phi_1 : X \to A$ is a homotopy inverse of the inclusion $A \subset X$. 

\begin{lemma}\label{lem_HD}
{\rm (1)} Suppose a subspace $A$ is HD in $X$. Then 
\begin{itemize}
\item[] \begin{itemize}
\item[{\rm (i)}\,] the inclusion $A \subset X$ is a $($controlled$)$ homotopy equivalence, and 
\item[{\rm (ii)}] $X$ is an ANR iff $A$ is an ANR. 
\end{itemize}
\item[{\rm (2)}] A subspace $A$ is HD in $X$ iff 
every $x \in X$ admits an open neighborhood $U$ in $X$ and 
a homotopy $\phi_t : U \to X$ such that $\phi_0$ 
is the inclusion $U \subset X$ and $\phi_t(U) \subset A$ $(0 < t \leq 1)$. $($cf. \cite[Fact 4.1 (i)]{Ya1}$)$. 
\end{itemize}
\end{lemma}

\noindent Lemma~\ref{lem_HD}\,(1)(ii) follows from Lemma~\ref{lem_ANR}\,(5). 

We conclude this preliminary section with the following characterization of $\ell_2$-manifold topological groups \cite{DT, To}. 
An $\ell_2$-manifold is a separable metrizable space locally homeomorphic to $\ell_2$. 

\begin{theorem}\label{thm_top-gp}  
{\rm (1)} A topological group is an $\ell_2$-manifold iff it is a separable, non locally compact, 
\begin{itemize}
\item[] completely metrizable ANR. 

\item[{\rm (2)}] Two $\ell_2$-manifolds are homeomorphic iff they are 
homotopy equivalent. 
\end{itemize}
\end{theorem}

%%%%%%%%%%%%%%%%%%%%%%%

\section{Transformation groups with weak topology} 

This section includes some results on transformation groups 
with weak topology. In the next section, these results will be applied to the diffeomorphism groups of noncompact manifolds endowed with the compact-open $C^r$-topology. 

\subsection{Transformation groups} \mbox{} 

A {\em transformation group} 
means a pair $(G, M)$ in which 
$M$ is a locally compact, $\sigma$-compact Hausdorff space and 
$G$ is a topological group acting on $M$ continuously and effectively. 
Each $g \in G$ induces a homeomorphism of $M$,
 which is also denoted by the same symbol $g$. 
 Let $G^c = \{g \in H \mid {\rm supp}(g) \text{ is compact}\}$. 
For any subsets $K, N$ of $M$ we obtain the following subgroups of $G$:
\[ \mbox{$G_K = \{g \in G : g|_K = {\rm id}_K \}, \ \ G(N) = G_{M \setminus N}, 
\ \ G_K(N) = G_K \cap G(N), \ \ G_{K}^c = G_K \cap G^c$ \ etc.}\]   
For any subgroup $H$ of $G$, 
let $H_0$ denote the connected component of the unit element $e$ in $H$. 
For example, the symbol $G^c_0$ denotes the connected component of $e$ in $G^c$. 
We also consider a subgroup $G^{c\,\ast}_0$ of $G^c_0$ defined by 
\[ \mbox{$G^{c\,\ast}_0 = \big\{ h \in G^c \mid h \in G(K)_0$ for some compact subset $K$ of $M$ $\big\}$.} \] 

For subsets $K \subset L \subset N$ of $M$, 
we have the set of embeddings 
$${\cal E}_K^G(L, N) = \big\{ g|_L : L \to M \mid g \in G_K(N) \big\}.$$
(If $K = \emptyset$, the symbol $K$ is omitted from the notation.) 
The group $G_K(N)$ acts transitively on the set ${\cal E}_K^G(L, N)$ 
by $g \cdot f = gf$ ($g \in G_K(N)$, $f \in {\cal E}_K^G(L, N)$). 
The restriction map 
\[ \mbox{$r : G_K(N) \to {\cal E}_K^G(L, N)$, \quad $r(g) = g|_L.$} \] 
coincides with the orbit map at the inclusion $i_L : L \subset M$ under this action. 

We need to pay an attention on the topology of the space ${\cal E}_K^G(L, N)$. 
A topology on ${\cal E}_K^G(L, N)$ is called {\em admissible} if 
the $G_K(N)$-action is continuous with respect to this topology. 
The strongest admissible topology on ${\cal E}_K^G(L, N)$ is the quotient topology induced by the map $r$. 
Otherwise specified, the set ${\cal E}_K^G(L, N)$ is endowed with this quotient topology. 
For $K \subset L_1 \subset L \subset N_1 \subset N$, one sees that ? 
${\cal E}_K^G(L, N_1)$ is a subspace of ${\cal E}_K^G(L, N)$ and 
the restriction map ${\cal E}_K^G(L, N) \to {\cal E}_K^G(L_1, N)$ is continuous. 

When the set ${\cal E}_K^G(L, N)$ is endowed with a specific admissible topology $\tau$, we write ${\cal E}_K^G(L, N)^\tau$ to avoid the ambiguity. 
The map $r : G_K(N) \to {\cal E}_K^G(L, N)^\tau$ is continuous. 
For any subset ${\cal F}$ of ${\cal E}_K^G(L, N)^\tau$, let 
${\cal F}^\tau$ denote the space ${\cal F}$ endowed with the subspace topology induced from $\tau$. When $i_L \in {\cal F}$, 
let ${\cal F}_0^\tau$ denote the connected component of $i_L$ in ${\cal F}^\tau$. 

We say that a map $f : X \to Y$ has a local section at $y \in Y$ if there exists a neighborhood $U$ of $y$ in $Y$ and a map $s : U \to X$ with $fs = i_U$. 

\begin{lemma}\label{lem_loc-sec}
Suppose $H$ is a subgroup of $G_K(N)$ and the restriction map $r : G_K(N) \to {\cal E}_K^G(L, N)^\tau$ 
has a local section $s : {\cal U} \to H \subset G_K(N)$ at $i_L$. 
Then, the following hold. 
\begin{itemize}
\item[{\rm (1)}] The restriction map $r|_H : H \to {\cal E}^H(L, N)^\tau$ is a principal bundle with the structure group $H_L$. 
\item[{\rm (2)}] 
\begin{itemize}
\item[(i)\ ] The topology $\tau$ coincides with the quotient topology. 
\item[(ii)\,] ${\cal E}^H(L, N)^\tau$ is open in ${\cal E}_K^G(L, N)^\tau$. 
\item[(iii)] If $H$ is a normal subgroup of $G_K(N)$, then  
\begin{itemize}
\item[(a)] each orbit of $H$ is closed and open in ${\cal E}_K^G(L, N)^\tau$ \\
$($in particular, ${\cal E}^H(L, N)^\tau$ is closed and open in ${\cal E}_K^G(L, N)^\tau)$, 
\item[(b)] if $H \subset G_K(N)_0$, then ${\cal E}^H(L, N)^\tau = {\cal E}_K^G(L, N)^\tau_0$. 
\end{itemize}
\end{itemize}
\end{itemize}
\end{lemma} 

\begin{proof}
(1) Note that (a) ${\cal E}^H(L, N)^\tau = H \cdot i_L$ and the map $r|_H$ coincides with the orbit map at $i_L$ under the action of the group $H$ and 
(b) ${\cal U}\subset {\cal E}^H(L, N)^\tau$ and $s$ is a local section of $r|_H$ at $i_L$. Hence, the statement (1) follows from the well known fact on the orbit map in the theory of group action.  

(2)\,(i) By (1) the map $r$ itself is a principal bundle. 
Hence, the map $r$ is a quotient map and the topology $\tau$ coincides with the quotient topology. 

(ii) We may assume that ${\cal U}$ is open in ${\cal E}_K^G(L, N)^\tau$. 
Since ${\cal U} \subset {\cal E}^H(L, N)^\tau$, we have ${\cal E}^H(L, N)^\tau = H {\cal U}$, which is open in ${\cal E}_K^G(L, N)^\tau$. 

(iii) Since $H i_L$ is open, the orbit $H(g \cdot i_L) = g(H i_L)$ is also open in ${\cal E}_K^G(L, N)^\tau$. 
\end{proof} 

In many cases, the set ${\cal E}_K^G(L, N)$ admits a natural admissible topology $\tau$ (for instance, the compact-open $C^r$-topology $(r=0, 1, \cdots, \infty)$). 
However, this topology $\tau$ coincides with the quotient topology, 
once we obtain a bundle theorem under the topology $\tau$.
Hence, our convention does not lose a generality of the arguments in the subsections  below. 

Finally we extract a behavior of the compact-open $C^\infty$ topology on diffeomorphism groups of a noncompact manifold $M$ at the ends of $M$ 
extend the notion of the compact-open $C^\infty$ topology on diffeomorphism groups to transformation groups. 

\begin{defn}~\label{weak_top} We say that a transformation group $(G, M)$ has a {\em weak} topology if it satisfies the following condition: 
\begin{itemize} 
\item[$(\ast)$] For any neighborhood $U$ of $e$ in $G$ there exists a compact subset $K$ of $M$ such that $G_K \subset U$.  
\end{itemize} 
\end{defn} 

\begin{remark} (1) 
If $G$ admits a compatible metric $\rho$, then $(G, M)$ has a weak topology iff 
for any $\varepsilon > 0$ there exists a compact subset $K$ of $M$ such that 
${\rm diam}_\rho \, G_K < \varepsilon$. We can always take $\rho$ to be left-invariant. 

(2) If a transformation group $(G, M)$ has a weak topology, then so is $(H, M)$ for any subgroup $H$ of $G$.  
\end{remark} 

The notion of weak topology is an extension of the notion of 
compact-open topology to transformation groups. 
It is readily seen that 
the compact-open $C^\infty$ topology on the diffeomorphism groups of smooth manifolds is also an example of weak topology. 

%%%%%%%%%%%%%%%%%%%%%%

\subsection{Basic assumptions on transformation groups} \mbox{} 

Suppose $(G, M)$ is a transformation group. 
Since $M$ is locally compact and $\sigma$-compact, 
there exists a sequence $\{ M_i \}_{i \geq 1}$ of compact subsets of $M$ such that 
\[ \mbox{$M = \cup_{i=0}^{\infty} \, M_i$ \ \ and \ \ $M_{i-1} \subset {\rm Int}_M M_i$ ($i \geq 1$), \ where $M_0 = \emptyset$.} \] 
This sequence is called an {\em exhausting sequence} of $M$.
Let $U_i = {\rm Int}_M\,M_i$ ($i \geq 1$). 
Consider the following basic conditions on the tuple $(G, M, \{ M_i \}_{i \geq 1})$.

\begin{assumption}\label{assumption} \mbox{} 
\begin{itemize}
\item[{\rm (A-0)}] The group $G$ is metrizable and the transformation group $(G, M)$ has a weak topology. 
\vskip 1mm 
\item[{\rm (A-1)}] 
For each $j > i > k \geq 0$, the restriction map 
$$\pi^i_{k,j} : G_{M_k}(U_j)_0 \longrightarrow {\cal E}^G_{M_k}(M_i, U_j)_0$$ 
is a principal bundle with the structure group \ \ 
$\displaystyle {\cal G}^i_{k,j} \equiv G_{M_k}(U_j)_0 \cap G_{M_i}$. 
\item[{\rm (A-2)}] 
\begin{itemize}
\item[(i)\,] $G(U_i)_0$ is an ANR for each $i \geq 1$.  
\vskip 1mm  
\item[(ii)] The spaces ${\cal U}_{k,j}^i = {\cal E}^G_{M_k}(M_i, U_j)_0$ $(j > i > k \geq 0)$ satisfy the next conditions: \\[1mm]  
${\cal U}_{k,j}^i$ is an open subspace of ${\cal E}^G_{M_k}(M_i, M)_0$, \ 
${\cal E}^G_{M_k}(M_i, M)_0 = \cup_{j > i}\,{\cal U}_{k,j}^i$ \ and \ $cl\,{\cal U}_{k,j}^i \subset {\cal U}_{k,j+1}^i$.
\end{itemize}
\end{itemize} 
\end{assumption} 

Below we assume that the tuple $(M, G, \{ M_i\}_{i \geq 1})$ satisfies the assumption (A). 

\begin{lemma}\label{lem_(A-3)} 
{\rm (1)} For each $i > k \geq 0$, the restriction map 
$$\pi^i_k : (G_{M_k})_0  \longrightarrow {\cal E}^G_{M_k}(M_i, M)_0$$ 
\begin{itemize}
\item[] is a principal bundle with the structure group \ \ 
$\displaystyle {\cal G}^i_k \equiv (G_{M_k})_0 \cap G_{M_i}.$

\item[(2)] The spaces $G_{M_k}(U_i)_0$ and ${\cal E}^G_{M_k}(M_i, M)_0$ are 
ANR's for each $i > k \geq 0$. 
\end{itemize} 
\end{lemma}

\begin{proof}
(1) This lemma follows from (A-1), (A-2)(ii) and Lemma~\ref{lem_loc-sec}.

(2) (i) By (A-1), for each $i > k \geq 1$, the restriction map 
$$\pi^k_{0,i} : G(U_i)_0 \longrightarrow {\cal E}^G(M_k, U_i)_0$$ 
is a principal bundle with the fiber \ 
${\cal G}^k_{0,i} \equiv G(U_i)_0 \cap G_{M_k}$. 
Since $G(U_i)_0$ is an ANR by (A-2)(i), so are ${\cal G}^k_{0,i}$ and 
$G_{M_k}(U_i)_0 = ({\cal G}^k_{0,i})_0$. 

(ii) By (i) and (A-1) ${\cal E}^G_{M_k}(M_i, U_j)_0$ is an ANR for each $j > i$,  
and so is ${\cal E}^G_{M_k}(M_i, M)_0$ by (A-2)(ii). 
\end{proof} 

\begin{lemma}\label{lem_submfd} {\rm (1)}{\rm (i)} 
 If $G_0$ is an ANR, then so is $(G_{M_i})_0$ for any $i \geq 1$. 
\begin{itemize} 
\item[] {\rm (ii)} If $G_{M_i}$ is an ANR for some $i \geq 1$, then so is $G_0$.  
\item[(2)] If $(G_{M_i})_0^{c \, \ast}$ is HD in $(G_{M_i})_0$ and 
$(G_{M_i})_0$ is open in $G_{M_i}$ for some $i \geq 1$, then 
${G_0^c}^\ast$ is HD in $G_0$. 
\end{itemize}
\end{lemma} 

\begin{proof} (1) 
The restriction map \ 
$\pi_0^i : G_0 \to {\cal E}^G(M_i, M)_0$ \ 
is a bundle map with fiber ${\cal G}_0^i = G_0 \cap G_{M_i}$. 
Hence, there exists an open neighborhood ${\cal U}$ of 
the inclusion map $i_{M_i}$ in ${\cal E}^G(M_i, M)_0$ such that 
$(\pi_0^i)^{-1}({\cal U}) \cong {\cal U} \times {\cal G}_0^i$. 
Since ${\cal U}$ is an ANR by (A-1)(i), it follows that $(\pi_0^i)^{-1}({\cal U})$ is an ANR if and only if ${\cal G}_0^i$ is an ANR. The assertions (i) and (ii) follow from these observations. 

(2) By Lemma~\ref{lem_HD}\,(2) it suffices to verify the following assertion: 
\begin{itemize}
\item[$(\#)$] Every $h \in G_0$ admits an open neighborhood ${\cal V}$ in $G_0$ and a homotopy $\phi : {\cal V} \times [0, 1] \to G_0$ such that $\phi_0$ 
is the inclusion ${\cal V} \subset G_0$ and $\phi_t({\cal V}) \subset {G^c_0}^\ast$ $(0 < t \leq 1)$.  
\end{itemize}

Under the map $\pi^i_0 : G_0 \to {\cal E}^G(M_i, M)_0$, 
each $h \in G_0$ induces $h|_{M_i} \in {\cal E}^G(M_i, M)_0$, and by (A-1)(ii) we have $h|_{M_i} \in {\cal E}^G(M_i, U_j)_0$ for some $j > i$. 
By (A-2) the restriction map 
$$\pi^i_{0,j} : G(U_j)_0 \longrightarrow {\cal E}^G(M_i, U_j)_0$$ 
is a bundle map and so 
there exists an open neighborhood ${\cal U}$ of $h|_N$ in ${\cal E}^G(M_i, U_j)_0$ 
and a local section $s : {\cal U} \to G(U_j)_0$ of $\pi^i_{0,j}$ such that $s(h|_{M_i}) = h$. 
Choose a small open neighborhood ${\cal V}$ of $h$ in $G_0$ such that 
$\pi^i_0({\cal V}) \subset {\cal U}$. 
For each $g \in {\cal V}$ we have $s(g|_{M_i})|_{M_i} = g|_{M_i}$ and $s(g|_{M_i})^{-1} g \in G_{M_i}$. 
Since $s(h|_{M_i})^{-1} h = id_M \in (G_{M_i})_0$ and 
$(G_{M_i})_0$ is open in $G_{M_i}$, 
by replacing ${\cal V}$ by a smaller one, we may assume that 
$s(g|_{M_i})^{-1} g \in (G_{M_i})_0$ $(g \in {\cal V})$. 

There exists an absorbing homotopy $\psi_t$ of $(G_{M_i})_0$ into ${(G_{M_i})^c_0}^\ast$. 
Then the required homotopy $\phi_t : {\cal V} \to G_0$ is defined by  
\[ \mbox{$\phi_t(g) = s(g|_{M_i}) \psi_t(s(g|_{M_i})^{-1} g)$.} \] 
Then $\phi_0(g) = g$ and 
for $0 < t \leq 1$ we have $\phi_t(g) \in {G^c_0}^\ast$ since 
$s(g|_{M_i}) \in G(U_j)_0 \subset {G^c_0}^\ast$ and 
$\psi_t(s(g|_{M_i})^{-1} g) \in {(G_{M_i})^c_0}^\ast \subset {G^c_0}^\ast$. 
\end{proof} 

%%%%%%%%%%%%%%%%%%%%

\subsection{Contractibility conditions} \mbox{} 

In this subsection we deduce some conclusions under some contractibility conditions. 
Consider the following conditions on the tuple $(M, G, \{ M_i\}_{i \geq 1})$: 

\begin{condition} \mbox{} 
\begin{itemize}
\item[{(C-1)}] ${\cal E}^G_{M_k}(M_i, M)_0 \simeq \ast$ for each $i > k \geq 0$. 

\vskip 1mm
\item[{(C-2)}] ${\cal G}^i_{0,j} \simeq \ast$ for each $j > i \geq 1$. 
\end{itemize} 
\end{condition}

Below we assume that the tuple $(G, M, \{ M_i\}_{i \geq 1})$ satisfies the assumption (A). 
Since $G$ is metrizable by the condition (A-0), it admits a left-invariant metric $\rho$.

\begin{lemma}\label{lem_(C-1)} 
If $(G, M, \{ M_i\}_{i \geq 1})$ satisfies the condition {\rm (C-1)}, 
then the following hold. 
\begin{itemize} 
\item[{\rm (1)}] For each $i > k \geq 0$, 
\begin{itemize}
\item[(i)\,] the bundle $\pi^i_k$ is a trivial bundle and ${\cal G}_k^i = (G_{M_i})_0$,  
\item[(ii)] $(G_{M_k})_0$ strongly deformation retracts onto $(G_{M_i})_0$. 
\end{itemize}
\vskip 1mm
\item[{\rm (2)}] $(G_{M_k})_0$ is an AR for each $k \geq 0$. 
\end{itemize} 
\end{lemma} 

\begin{proof} 
The conditions (A-1) and (C-1) imply that ${\cal E}^G_{M_k}(M_i, M)_0$ 
is an AR. 

(1) Since the base space ${\cal E}^G_{M_k}(M_i, M)_0$ is metrizable and contractible, the bundle $\pi^i_k$ is trivial. 
This means that there exists a fiber-preserving homeomorphism over ${\cal E}^G_{M_k}(M_i, M)_0$, 
$$(G_{M_k})_0  \cong {\cal E}^G_{M_k}(M_i, M)_0 \times {\cal G}^i_k.$$
Since $(G_{M_k})_0$ is connected, so is ${\cal G}_k^i$ and we have ${\cal G}_k^i = (G_{M_i})_0$. Since the base ${\cal E}^G_{M_k}(M_i, M)_0$ is an AR, it admits a strong deformation retraction (SDR) onto the singleton $\{ i_{M_i} \}$. This induces the required SDR of $(G_{M_k})_0$ onto ${\cal G}_k^i$
\vskip 1mm 
(2)(i) First we show that $(G_{M_k})_0 \simeq \ast$. 
By (1)(ii), for each $i \geq 0$ 
there exists a SDR $h^i_t$ ($0 \leq t \leq 1$) of $(G_{M_i})_0$ onto $(G_{M_{i+1}})_0$. 
A SDR $h_t$ ($k \leq t \leq \infty$) of $(G_{M_k})_0$ onto $\{ id_M \}$ is defined by
\begin{eqnarray*}
 h_t(f) \ \, &=& h_{t-i}^ih_1^{i-1} \cdots h_1^k(f) \ \ \ (f \in (G_{M_k})_0, \ i \geq k, \ i \leq t \leq i+1) \\
 h_\infty(f) &=& id_M.
\end{eqnarray*}
Since $(G, M)$ has a weak topology, it follows that 
${\rm diam}_\rho \, (G_{M_i})_0 \to 0$ and the homotopy $h : (G_{M_k})_0 \times [k, \infty] \to (G_{M_k})_0$ is continuous. 

(ii) To see that $(G_{M_k})_0$ is an ANR, 
we apply Lemma~\ref{lem_ANR}\,(5) (the Hanner's criterion). 
By (1), for each $i > k$ the restriction map 
$$\pi^i_k : (G_{M_k})_0 \longrightarrow {\cal E}^G_{M_k}(M_i, M)_0$$  
is a trivial bundle with an ANR base space and the fiber $\displaystyle {\cal G}^i_k \equiv (G_{M_i})_0$. 
By (i) the fiber ${\cal G}^i_k$ is contractible. 
These imply that $\pi^i_k$ admits a section $s^i_k$ and that $s^i_k \pi^i_k$ is $\pi^i_k$-fiber preserving homotopic to $id$. 
Since the metric $\rho$ is left-invariant, the diameter of each fiber of $\pi^i_k$ coincides with ${\rm diam}_\rho G_{M_i}$. 
Since ${\rm diam}_\rho G_{M_i} \to 0$ ($i \to 0$), 
the Hanner's criterion implies that $(G_{M_k})_0$ is an ANR. 
\end{proof}

\begin{lemma}\label{lem_(C-2)}
If $(G, M, \{ M_i\}_{i \geq 1})$ satisfies the condition {\rm (C-2)}, then 
${G^c_0}^\ast$ is HD in $G_0$. 
\end{lemma} 

\begin{proof} 
(1) For notational simplicity, for each $j > i \geq 1$ we set 
${\cal D}_j = G(U_j)_0$ and ${\cal U}_{i,j} = {\cal E}^G(M_i, U_j)_0$. 
By (A-2) the map 
$\pi^i_{0,j} : {\cal D}_j \to {\cal U}_{i,j}$ is a principal bundle. 
Since the fiber is an AR by (C-2) and the base space is metrizable, 
this bundle has a global section. 
Thus, the map $\pi^i_{0,j}$ is a trivial bundle with an AR fiber and 
hence it has the following relative lifting property: 
\begin{itemize} 
\item[($\ast$)]
If $Y$ is a metric space, $B$ is a closed subset of $Y$ and  
$\phi : Y \to {{\cal U}_{i,j}}$ and $\phi_0 : B \to {\cal D}_j$ are maps with $\phi|_B = \pi^i_{0,j}\phi_0$, then 
there exists a map $\Phi : Y \to {\cal D}_j$ such that $\pi^i_{0,j}\Phi = \phi$ and $\Phi|_B = \phi_0$.
\end{itemize}

(2) Next consider the principal bundle
$$\pi^i_0 : G_0 \to {\cal E}^G(M_i, M)_0.$$
For each $j > i \geq 1$, 
we set ${{\cal V}_{i,j}} = (\pi^i_0)^{-1}({{\cal U}_{i,j}}) \subset G_0$. Then ${{\cal U}_{i,j}}$, ${{\cal V}_{i,j}}$ and ${\cal D}_j$ satisfy the following conditions: 
for each $i \geq 1$ 
\begin{itemize}
\item[(i)\ ] ${\cal E}^G(M_i, M)_0 = \cup_{j > i} \, {\cal U}_{i,j}$, \ 
 ${\cal U}_{i,j}$ is open in ${\cal E}^G(M_i, M)_0$ \ and \ 
$cl \, {\cal U}_{i,j} \subset {\cal U}_{i,j+1}$. 
\item[(ii)\,] $G_0 = \cup_{j > i}\,{{\cal V}_{i,j}}$, \ 
 ${{\cal V}_{i,j}}$ is open in $G_0$, \ 
 $cl \, {{\cal V}_{i,j}} \subset {{\cal V}_{i,j+1}}$ \ and \  
 ${{\cal V}_{i+1,j}} \subset {{\cal V}_{i,j}}$ $(j > i+1)$. 
\item[(iii)]  ${G^c_0}^\ast = \cup_{j > i} \,{\cal D}_j$ \ and \ ${\cal D}_j \subset {\cal D}_{j+1}$. 
\end{itemize}

(3) We have to construct a homotopy $F : G_0 \times [0,1] \to G_0$ such that $F_0 = id$ and 
$F_t(G_0)$ $\subset {G^c_0}^\ast$ ($0 < t \leq 1$).
We replace the interval $[0,1]$ by $[1, \infty]$. 

(i) For each $i \geq 1$ we can find a map $s^i : {\cal E}^G(M_i, M)_0 \to {G^c_0}^\ast$ 
such that 
\[ \mbox{$s^i(f)|_{M_i} = f|_{M_i}$ ($f \in {\cal E}^G(M_i, M)_0$) and $s^i(cl \, {{\cal U}_{i,j}}) \subset {\cal D}_{j+1}$ ($j > i$).} \] 
In fact, using the property ($\ast$), 
inductively we can construct maps $s^i_j : cl \, {{\cal U}_{i,j}} \to {\cal D}_{j+1}$ ($j > i$) such that 
\[ \mbox{$s^i_j(f)|_{M_i} = f$ ($f \in cl \, {{\cal U}_{i,j}}$) and $s_{j+1}^i|_{cl \, {{\cal U}_{i,j}}} = s^i_j$.} \] 
The map $s^i$ is defined by $s^i|_{cl \, {{\cal U}_{i,j}}} = s^i_j$. 
Let $F_i = s^i \pi^i_0 : G_0 \to {G^c_0}^\ast$. 
We have $F_i(cl \, {{\cal V}_{i,j}}) \subset {\cal D}_{j+1}$ and $F_i(h)|_{M_i} = h|_{M_i}$. 

(ii) For each $i \geq 1$,    
we can inductively construct a sequence of homotopies 
\[ \mbox{$H^j : cl \, {{\cal V}_{i+1, j}} \times [i, i+1] \to {\cal D}_{j+1}$ ($j > i+1$)} \]   
such that 
\[ \mbox{$H^j_i = F_i$, $H^j_{i+1} = F_{i+1}$, $H^{j+1}|_{cl \, {{\cal V}_{i+1, j}} \times [i, i+1]} = H^j$ and $H^j_t(h)|_{M_i} = h|_{M_i}$.} \]  
If $H^j$ is given, then $H^{j+1}$ is obtained 
by applying the property ($\ast$) to the diagram: 
\vskip 3mm
\hspace*{10mm}
\begin{minipage}[c]{120pt}
\begin{tabular}[c]{ccc}
$B$ & 
$\stackrel{\phi_0}{\longrightarrow}$ & ${\cal D}_{j+2}$ \\[1mm]
$\bigcap$ & & \smash{\raisebox{-3pt}{$\big\downarrow$}} \\[1mm]
$Y$ & $\stackrel{\phi}{\longrightarrow}$ & ${\cal U}_{i, j+2}$, 
\end{tabular}
\end{minipage}
\begin{minipage}[c]{270pt}
$\phi(h,t) = h|_{M_i}$, \hspace{10pt}
$\phi_0(h, t) = 
\begin{cases} 
H^j(h,t) & (h \in cl \, {{\cal V}_{i+1, j}}) \\[1mm]
F_t(h) & (t = i, i+1) 
\end{cases}$
\end{minipage}
\vskip 1mm
\[ (Y, B) = (cl \, {{\cal V}_{i+1, j+1}} \times [i, i+1], (cl \, {{\cal V}_{i+1, j}} \times [i, i+1]) \cup (cl \, {{\cal V}_{i+1, j+1}} \times \{i, i+1 \})). \]
Thus we can define a homotopy $F : G_0 \times [i, i+1] \to {G^c_0}^\ast$ by $F = H^j$ on $cl \, {\cal V}_{i+1, j} \times [i, i+1]$. 
Since $F_t(h)|_{M_i} = h|_{M_i}$ for $t \geq i$, we can continuously extend $F$ by $F_{\infty} = id$. 
\end{proof} 

%%%%%%%%%%%%%%%%%%%

Lemmas~\ref{lem_(C-1)} and \ref{lem_(C-2)} 
yield the following criterions. 

\begin{theorem}\label{thm_criterion} 
Suppose the tuple $(G, M, \{ M_i\}_{i \geq 1})$ satisfies the assumption {\rm (A)}. 
\begin{itemize} 
\item[(1)] If $G_{M_k}(U_i)_0 \simeq \ast$ and 
${\cal G}^i_{k,j}$ is connected $($i.e., $G_{M_k}(U_j)_0 \cap G_{M_i} = G_{M_i}(U_j)_0$\,$)$ for each $j > i > k \geq 0$, then 
$G_0$ is an AR and ${G^c_0}^\ast$ is HD in $G_0$. 
\vskip 1mm 
\item[(2)] 
\begin{itemize}
\item[(i)\,]  If $G_{M_i} \simeq \ast$ for each $i \geq 1$, then $G_0$ is an ANR.

\item[(ii)] If $G_{M_1}$ is connected and 
${\cal G}_{1,j}^i = G_{M_1}(U_j)_0 \cap G_{M_i} \simeq \ast$ for each $j > i \geq 2$, then 
$G^{c\,\ast}_0$ is HD in $G_0$. 
\end{itemize} 
\end{itemize} 
\end{theorem}

\begin{proof}
(1) By Lemmas~\ref{lem_(C-1)} and \ref{lem_(C-2)} it suffices to verify the following assertions: \ for each $j > i > k \geq 0$ 
\begin{itemize} 
\item[{\rm (i)}\,] the bundle $\pi^i_{k,j}$ is trivial, and the spaces 
${\cal G}^i_{k,j}$, $G_{M_k}(U_j)_0$ and ${\cal E}^G_{M_k}(M_i, U_j)_0$ are AR's, 
\vskip 1mm 
\item[(ii)]  
${\cal E}^G_{M_k}(M_i, M)_0$ is an AR. 
\end{itemize} 

(i) By Lemma~\ref{lem_(A-3)}\,(2) and the assumption, both the total space $G_{M_k}(U_j)_0$ and 
the fiber ${\cal G}^i_{k,j} = G_{M_i}(U_j)_0$ are AR's. 
From the homotopy exact sequence of 
the bundle $\pi^i_{k,j}$ it follows that the ANR base space 
${\cal E}^G_{M_k}(M_i, U_j)_0$ is contractible and the bundle $\pi^i_{k,j}$ is trivial. 
(Alternatively, since the fiber is an AR, the principal bundle $\pi^i_{k,j}$ admits a global section, which means that this bundle is trivial and the base space is an AR since it is a retract of the AR total space.) 

(ii) The assertion follows from (i), (A-2)(ii) and Lemma~\ref{lem_ANR}\,(3). 
\vskip 1mm 
(2) Since $(G_{M_1})_0 = G_{M_1}$, 
by Lemma~\ref{lem_submfd} 
it suffices to show the next assertions:  
\begin{itemize} 
\item[ (i$'$)\,] $(G_{M_1})_0$ is an ANR if $G_{M_i} \simeq \ast$ for each $i \geq 1$. 
\item[(ii$'$)] $(G_{M_1})_0^{c \, \ast}$ is HD in $(G_{M_1})_0$ if 
${\cal G}_{1,j}^i = G_{M_1}(U_j)_0 \cap G_{M_i} \simeq \ast$ for each $j > i \geq 2$. 
\end{itemize} 
Since the tuple $(G_{M_1}, M, \{ M_i\}_{i \geq 2})$ also 
satisfies the assumption (A), by Lemmas~\ref{lem_(C-1)} and \ref{lem_(C-2)} 
it remains to verify the condition (C-1) in (i$'$) and (C-2) in (ii$'$) respectively. 

(i$'$) ${\cal E}^G_{M_k}(M_i, M)_0 \simeq \ast$ for each $i > k \geq 1$: 
In fact, the restriction map 
$$\pi^i_k : (G_{M_k})_0  \longrightarrow {\cal E}^G_{M_k}(M_i, M)_0$$ 
is a principal bundle with the fiber \ 
${\cal G}^i_k \equiv (G_{M_k})_0 \cap G_{M_i}.$
Since $(G_{M_k})_0 = G_{M_k} \simeq \ast$ and 
${\cal G}^i_k = G_{M_i} \simeq \ast$, it follows that 
the ANR base space ${\cal E}^G_{M_k}(M_i, M)_0$ is also contractible. 

(ii$'$) ${\cal G}_{1,j}^i = G_{M_1}(U_j)_0 \cap G_{M_i} \simeq \ast$ for each $j > i \geq 2$: \ This is the condition (C-2) itself for the tuple $(G_{M_1}, M, \{ M_i\}_{i \geq 2})$. 
This completes the proof. 
\end{proof} 

%%%%%%%%%%%%%%%%%%%%%

\section{Diffeomorphism groups of non-compact $n$-manifolds}

\subsection{General properties of diffeomorphism groups of $n$-manifolds} \mbox{}

Suppose $M$ is a smooth (separable metrizable) $n$-manifold without boundary 
and $X$ is a closed subset of $M$ with $X \neq M$.
Let $r = 1, 2, \cdots, \infty$. 

When $N$ is a smooth manifold, the symbol ${\cal C}^r(M, N)$ denotes the space of $C^r$-maps $f : M \to N$ with the compact-open $C^r$-topology. 
For any map $f_0 : X \to N$, let ${\cal C}^{\, r}_{f_0}(M, N)$ denote the subspace $\{ f \in {\cal C}^r(M, N) \mid f|_X = f_0 \}$. 
When $M$ is a smooth submanifold of $N$, 
the symbol ${\cal E}^{\, r}_X(M, N)$ denotes the subspace of ${\cal C}^r(M, N)$ consisting of $C^{\, r}$-embeddings $f : M \hookrightarrow N$ with $f|_X = id_X$ 
and ${\cal E}^{\, r}_X(M, N)_0$ denotes the connected component of the inclusion $i_N : N \subset M$ in ${\cal E}^{\, r}_X(N, M)$. 
For spaces $Y$ and $Z$, the symbol ${\cal C}^0(Y, Z)$ denotes the space of $C^0$-maps $f : Y \to Z$ with the compact-open topology. 

Consider the jet-map $j^r : C^r(M, N) \to C^0(M, J^r(M, N))$. 
It is a closed embedding \cite[Ch2. Section 4, p.61--62]{Hi}. 
Thus, if we choose a complete metric $d$ on the jet-bundle $J^r(M, N)$ with $d \leq 1$ and 
a sequence of compact $n$-submanifolds $M_i$ ($i \geq 1$) such that $M_i \subset {\rm Int}\,M_{i+1}$ and $M = \cup_i\,M_i$, 
then we can define a complete metric $\rho$ on ${\cal C}^r(M, N)$ by 
\[ \rho(f, g) = \sum_{i=1}^{\infty} \, \frac{1}{2^i} \sup_{x \in M_i} \, d(j^{\, r}_x f, j^{\, r}_x g). \]
When $N = M$, it induces a metric on ${\cal D}^{\, r}_X(M) \subset {\cal C}^r(M, M)$. 
We can define a complete metric $\rho^\ast$ on ${\cal D}^{\, r}_X(M)$ by  
\[ \rho^\ast(f, g) = \rho(f, g) + \rho(f^{-1}, g^{-1}). \]

\begin{lemma}\label{lem_diff-gp} 
{\rm (1)} ${\cal D}^{\, r}_X(M)$ is a topological group \cite[Ch\,2.,Section 4.p.~64]{Hi}. 
\begin{itemize}
\item[{\rm (2)}] ${\cal D}^{\, r}_X(M)$ and ${\cal D}^{\, r}_X(M)_0$ are separable, completely metrizable and not locally compact 
\cite[Ch\,2., Section 4, p.~61--62, Theorems 4.3, 4.4]{Hi}. 
\vskip 1mm
\item[(3)] ${\cal D}^{\, r}_X(M)$ and ${\cal D}^{\, r}_X(M)_0$ are $\ell^2$-manifolds iff they are ANR's.  
\end{itemize}
\end{lemma}

If $M$ is a closed smooth $n$-manifold, then ${\cal D}^r(M)$ is a smooth Fr\'echet manifold \cite[Section I.4, Example 4.1.3, etc.]{Ham}.
In this paper we are only concerned with topological Fr\'echet manifolds (= topological $\ell_2$-manifolds). 
Below, the emphasis is put on relative cases. 

Suppose $M$ is a compact smooth $n$-manifold (possibly with boundary) and $X$ is a closed subset of $M$. 
When $E \to M$ is a smooth fiber bundle over $M$ and $s_0 : X \to E$ is a section of $E$ over $X$, 
the symbol $\Gamma^{\, r}_{s_0}(M, E)$ denotes the space of $C^r$-sections $s$ of $E$ over $M$ 
with $s|_X = s_0$, endowed with the compact-open $C^r$-topology. 
When $E$ is a vector bundle over $M$, the space $\Gamma^{\, r}_{s_0}(M, E)$ is a Fr\'echet space.
If $E$ is a fiber bundle over $M$ and the fiber of $E$ is a smooth manifold without boundary, then $\Gamma^{\, r}_{s_0}(M, E)$ is a Fr\'echet manifold.
A local chart around $s \in \Gamma^{\, r}_{s_0}(M, E)$ is defined as follows: 
Consider the restriction $T^v(E)|_{s(M)}$ of the vertical tangent bundle $T^v(E)$ of $E$ to the image $s(M)$ in $E$. 
There exists a diffeomorphism between an open neighborhood of $s(M)$ in $E$ and an open neighborhood of the image of zero section $0$ in $T^v(E)|_{s(M)}$. 
It induces a homeomorphism between an open neighborhood of $s$ in $\Gamma^{\, r}_{s_0}(M, E)$ and 
an open neighborhood of $0$ in $\Gamma^{\, r}_{0_X}(s(M), T^v(E)|_{s(M)})$.
If $N$ is a smooth manifold without boundary and $f_0 : X \to N$ is a map, 
then the space ${\cal C}^{\, r}_{f_0}(M, N) = \{ f \in {\cal C}^r(M, N) \mid f|_X = f_0 \}$ 
is also a Fr\'echet manifold since it is identified with $\Gamma^{\, r}_{s_0}(M, E)$ for the trivial bundle $E = M \times N \to M$ and the section $s_0(x) = (x, f_0(x))$ over $X$

\begin{lemma}\label{lem_Frechet} 
{\rm (i)} Suppose $N$ is a smooth manifold without boundary, $M$ is a comapct smooth submanifold of $N$ and $X$ is a closed subset of $M$. 
Then ${\cal E}^{\, r}_X(M, N)$ is a Fr\'echet manifold. \\
{\rm (ii)} Suppose $M$ is a compact smooth $n$-manifold and $X$ is a closed subset of $M$. 
If $\partial M \subset X$ or $\partial M \cap X = \emptyset$, then ${\cal D}^{\, r}_X(M)$ is a Fr\'echet manifold. 
\end{lemma}

\begin{proof}
(i) ${\cal E}^{\, r}_X(M, N)$ is open in ${\cal C}^{\, r}_{i_X}(M, N)$ \cite[Ch.2 Theorem 1.4]{Hi}. 

(ii) If $\partial M \subset X$, then ${\cal D}^{\, r}_X(M) = {\cal E}^{\, r}_X(M, M) = {\cal E}^{\, r}_X(M, \tilde{M})$, 
where $\tilde{M}$ is the open manifold obtained from $M$ by attaching an open collar to $\partial M$. 
In the case $X \cap \partial M = \emptyset$, consider the restriction map $\pi : {\cal D}^{\, r}_X(M) \to {\cal D}^r(\partial M)$, $\pi(h) = h|_{\partial M}$. 
Using a collar of $\partial M$ in $M$, for any $f \in {\cal D}^r(\partial M)$ which is sufficiently close to $id$, 
we can easily construct a canonical extension $\tilde{f} \in {\cal D}^{\, r}_X(M)$ of $f$. 
This implies that ${\rm Im}\,\pi$ is clopen subset of ${\cal D}^r(\partial M)$ and 
${\cal D}^{\, r}_X(M) \to {\rm Im}\,\pi$ is a principal bundle with fiber ${\cal D}^{\, r}_{X \cup \partial M}(M)$. 
Since ${\cal D}^r(\partial M)$ and ${\cal D}^{\, r}_{X \cup \partial M}(M)$ are Fr\'echet manifolds, so is ${\cal D}^{\, r}_X(M)$. 
\end{proof}
 
\begin{remark}\label{rem_path}
(1) A family $h_t \in {\cal D}^{\, r}_X(M)$ ($t \in [0, 1]$) is called a $C^r$-isotopy rel $X$ 
if $H : M \times [0, 1] \to M$, $H(x, t) = h_t(x)$, is a $C^r$-map.
In this case $h_t$ is a path in ${\cal D}^{\, r}_X(M)$ (i.e., $[0, 1] \ni t \mapsto h_t \in {\cal D}^{\, r}_X(M)$ is continuous).
In general, if $h_t$ is a path in ${\cal D}^{\, r}_X(M)$, then $h_t$ is a $C^0$-isotopy rel $X$, but $H$ is not necessarily $C^1$ in $t$. 

(2) Since Fr\'echet manifolds are localy path-connected, 
the connected components ${\cal E}^{\, r}_X(N, M)_0$ and ${\cal D}^{\, r}_X(M)_0$ in Lemma 2.2 are path-connected. 
Thus, if $h \in {\cal D}^{\, r}_X(M)_0$, 
then there is a path $h_t$ ($t \in [0,1]$) in ${\cal D}^{\, r}_X(M)_0$ with $h_0 = id_M$ and $h_1 = h$, 
which is a $C^0$-isotopy rel $X$. 
\end{remark}

%%%%%%%%%%%%%%%%%%

\subsection{The bundle theorem} \mbox{}

The bundle theorem connecting 
diffeomorphism groups and embedding spaces \cite{Ce, KM, Pa1, Seeley} plays an essential role in order to apply Theorem~\ref{thm_criterion} (recall Assumption (A)).  
Suppose $M$ is a smooth $n$-manifold without boundary. 
A compact smooth submanifold of $M$ means 
a compact subset $N$ of $M$ which is the union of a disjoint family $\{ N_k \}_{k=0}^n$ such that $N_k$ is a (possibly empty) closed smooth $k$-submanifold of $M$ for $k = 0, 1, \cdots, n-1$ and $N_n$ is a (possibly empty) compact smooth $n$-submanifold of $M$. 

Suppose $N$ is a compact smooth submanifold of $M$ and $X$ is a closed subset of $N$. 
Let $U$ be any open neighborhood of $N$ in $M$.

\begin{theorem}\label{thm_extension}
For any $f \in {\cal E}^{\, r}_X(N, U)$ there exist 
a neighborhood ${\cal U}$ of $f$ in ${\cal E}^{\, r}_X(N, U)$ and a map $\phi : {\cal U} \to {\cal D}^{\, r}_{X \cup (M \setminus U)}(M)_0$  
such that $\phi(g)f = g$ $(g \in {\cal U})$ and $\phi(f) = id_M$. 
\end{theorem}

Consider the restriction map $\pi : {\cal D}^{\, r}_{X \cup (M \setminus U)}(M) \to {\cal E}^{\, r}_X(N, U)$, $\pi(h) = h|_N$. 
The group ${\cal D}^{\, r}_{N \cup (M \setminus U)}(M)$ acts on ${\cal D}^{\, r}_{X \cup (M \setminus U)}(M)$ by right composition. 

\begin{corollary}\label{cor_bdle}
{\rm (1)} The image $\pi({\cal D}^{\, r}_{X \cup (M \setminus U)}(M))$ is open and closed in ${\cal E}^{\, r}_X(N, U)$ and the map 
$$\pi : {\cal D}^{\, r}_{X \cup (M \setminus U)}(M) \to \pi({\cal D}^{\, r}_{X \cup (M \setminus U)}(M))$$ 
\begin{itemize}
\item[] is a principal bundle with fiber ${\cal D}^{\, r}_{N \cup (M \setminus U)}(M)$. 
\item[{\rm (2)}] The restriction map 
$$\pi : {\cal D}^{\, r}_{X \cup (M \setminus U)}(M)_0 \to {\cal E}^{\, r}_X(N, U)_0, \ \pi(h) = h|_N,$$ 
is a principal bundle with fiber ${\cal D}^{\, r}_{X \cup (M \setminus U)}(M)_0 \cap {\mathcal D}_N(M)$.  
\end{itemize}
\end{corollary}

\subsection{Diffeomorphism groups of noncompact $n$-manifolds} \mbox{} 

Suppose $M$ is a noncompact connected smooth $n$-manifold without boundary and $X$ is a compact smooth submanifold of $M$. 
A {\em smooth} exhausting sequence of $(M, X)$ means an exhausting sequence 
$\{ M_i \}_{i \geq 1}$ of $M$ such that each $M_i$ is a compact $n$-submanifold of $M$ for each $i \geq 1$ and $X \subset {\rm Int}\,M_1$. Let $M_0 = X$. 
Obviously $(M, X)$ has a smooth exhausting sequence and 
any smooth exhausting sequence of $(M, X)$ satisfies the assumption (A) 
with respect to the diffeomorphism group ${\cal D}_X(M)$. 
Therefore, Theorem~\ref{thm_criterion} and Lemma~\ref{lem_diff-gp}\,(3) yield the following consequences. 

\begin{proposition}\label{prop_criterion_diff} Suppose $M$ admits a smooth exhausting sequence $\{ M_i \}_{i \geq 1}$ such that 
\begin{itemize}
\item[(i)\,] ${\cal D}_{M_k\cup (M - U_i)}(M)_0 \simeq \ast$ \ for each $i > k \geq 0$, 
\item[(ii)] ${\cal D}_{M_k\cup (M - U_j)}(M)_0 \cap {\cal D}_{M_i}(M) 
= {\cal D}_{M_i\cup (M - U_j)}(M)_0$ \ for each $j > i > k \geq 0$. 
\end{itemize}
Then, {\rm (1)} ${\cal D}_X(M)_0 \cong \ell^2$ \ and \ {\rm (2)} 
${\cal D}_X(M)^{c\,\ast}_0$ is HD in ${\cal D}_X(M)_0$. 
\end{proposition} 

\begin{proposition}\label{prop_criterion_cpt} If $M = {\rm Int}\,N$ for some compact $n$-manifold $N$ with non-empty boundary, then 
\begin{itemize}
\item[(1)] ${\cal D}_X(M)_0$ is an $\ell^2$-manifold,  
\item[(2)] if ${\cal D}_{\partial N \times \{ 0 \}}(\partial N \times [0,1]) \simeq \ast$, then ${\cal D}_X(M)_0 \cong {\cal D}_X(N)_0$, 
\item[(3)] if ${\cal D}_{\partial N \times \{ 0, 1 \}}(\partial N \times [0,1])_0 \simeq \ast$, then ${\cal D}_X(M)^{c\,\ast}_0$ is HD in ${\cal D}_X(M)_0$.
\end{itemize} 
\end{proposition} 

\begin{proof} 
Take a smooth closed collar $\partial N \times [0,2]$ of $\partial N = \partial N \times \{ 0 \}$ in $N - X$ and let $M_i = M - \big(\partial N \times [0, 1/i)\big)$ $(i \geq 1)$. 
Then $\{ M_i \}_{i \geq 1}$ is a smooth exhausting sequence of $(M, X)$.
Let $G = {\cal D}^{\, r}_X(M)$. 
The tuple $(G, M, \{ M_i \}_{i \geq 1})$ satisfies the assumption (A). 

(1) By Theorem~\ref{thm_criterion}\,(2)(i) and Lemma~\ref{lem_diff-gp}\,(3) it suffices to show that $G_{M_i} = {\cal D}^{\, r}_{M_i}(M) \simeq \ast$ for each $i \geq 1$. 
This is verified by the Alexander trick towards $\infty$, since 
$M - {\rm Int}\,M_i = \partial N \times (0, 1/i]$ is an open collar of $\partial M_i = \partial N \times \{ 1/i \}$. 

(2) By (1) and Lemma~\ref{lem_Frechet}, both ${\cal D}^{\, r}_X(M)_0$ and ${\cal D}^{\, r}_X(N)_0$ are $\ell^2$-manifolds. 
Thus, by Theorem~\ref{thm_top-gp}\,(2) it suffices to show that ${\cal D}^{\, r}_X(M)_0 \simeq {\cal D}^{\, r}_X(N)_0$. 

Note that (i) $G_{M_1}$ is an AR and $(G_{M_1})_0 = G_{M_1}$ by (1) and (ii) 
${\cal D}^{\, r}_{M_1}(N)$ is an AR since it is an ANR and  
${\cal D}^{\, r}_{M_1}(N) \simeq {\cal D}_{\partial N \times \{ 1 \}}(\partial N \times [0,1]) \simeq \ast$ by the assumption. 
The restriction maps 
\[ \mbox{$\pi : {\cal D}^{\, r}_X(M)_0 \longrightarrow {\cal E}^{\, r}_X(M_1, M)_0$ \ \ and \ \ 
$\pi_1 : {\cal D}^{\, r}_X(N)_0 \longrightarrow {\cal E}^{\, r}_X(M_1, M)_0$}\] 
are principal bundles with the AR fibers 
\[ \mbox{$G_0 \cap G_{M_1} = G_{M_1}$ \ \ and \ \ 
${\cal D}^{\, r}_X(N)_0 \cap {\cal D}^{\, r}_{M_1}(N) = {\cal D}^{\, r}_{M_1}(N)$
.} \] 
Thus, by Lemma~\ref{lem_bdle_AR} the maps $\pi$ and $\pi_1$ are homotopy equivalences. 

(3) Let ${\cal H} = {\cal D}^{\, r}(\partial N \times {\Bbb R})$. Then,  
for $j > i \geq 2$ we have 
\begin{align*}
G_{M_1}(U_j)_0 \cap G_{M_i} &\cong {\cal H}(\partial N \times [0,2])_0 \cap 
{\cal H}(\partial N \times [1,2]) = {\cal H}(\partial N \times [1,2])_0 \\
&\simeq {\cal D}_{\partial N \times \{ 0, 1 \}}(\partial N \times [0,1])_0 \simeq \ast.
\end{align*}
Since $(G_{M_1})_0 = G_{M_1}$,  
the assertion follows from Theorem~\ref{thm_criterion}\,(2)(ii). 
\end{proof} 

Proposition~\ref{prop_finite-type} now follows from Proposition~\ref{prop_criterion_cpt}\,(1). 
In the next section, Propositions~\ref{prop_criterion_diff} and \ref{prop_criterion_cpt} are used to deduce 
Theorems~\ref{thm_l2} and \ref{thm_HD}. 

%%%%%%%%%%%%%%%%

\section{Diffeomorphism groups of non-compact $2$-manifolds} 

\subsection{Fundamental facts on diffeomorphism groups of $2$-manifolds} \mbox{} 

First we recall some fundamental facts on diffeomorphism groups of compact 2-manifolds. 
The symbols ${\Bbb S}^1$, ${\Bbb S}^2$, ${\Bbb T}$, ${\Bbb D}$, ${\Bbb A}$, ${\Bbb P}$, ${\Bbb K}$ and ${\Bbb M}$ denote 
the 1-sphere (circle), 2-sphere, torus, disk, annulus, projective plane, Klein bottle and M\"obius band, respectively. 

\begin{theorem}\label{thm_h-type_cpt}$($\cite{EE, Sm} etc.$)$ 
Suppose $M$ is a compact connected smooth $2$-manifold $($possibly with boundary$)$ and 
$X$ is a compact smooth submanifold of $M$. 
\begin{itemize}
\item[{\rm (i)}\ ] {\rm (a)} ${\cal D}^{\, r}(M)_0 \simeq SO(3)$ if $M \cong {\Bbb S}^2, {\Bbb P}$. \ 
{\rm (b)} ${\cal D}^{\, r}(M)_0 \simeq {\Bbb T}$ if $M \cong {\Bbb T}$. 
\item[{\rm (ii)}\,] 
${\cal D}^{\, r}_X(M)_0 \simeq {\Bbb S}^1$ if $(M,X) \cong  ({\Bbb S}^2, {\rm 1 pt}), ({\Bbb S}^2, {\rm 2 pt}), ({\Bbb D}, \emptyset), ({\Bbb D}, 0), ({\Bbb A}, \emptyset), ({\Bbb M}, \emptyset), ({\Bbb K}, \emptyset)$. 
\item[{\rm (iii)}] ${\cal D}^{\, r}_X(M)_0 \simeq \ast$ in all other cases. 
\item[{\rm (iv)}] ${\cal D}^{\, r}_\partial({\Bbb D}) \simeq {\cal D}^{\, r}_{\{ 0 \} \cup \partial}({\Bbb D}) \simeq \ast$  and ${\cal D}^{\, r}_\partial({\Bbb M}) \simeq \ast$. 
\end{itemize}
\end{theorem}

By \cite{EE} and a $C^r$-analogue of \cite{Ep} we have 

\begin{lemma}\label{lem_C^r-iso} Suppose $M$ is a compact smooth $2$-manifold $($possibly with boundary$)$. 
\begin{itemize}
\item[{\rm (1)}] Suppose $N$ is a closed collar of $\partial M$. 
If $h \in {\cal D}^{\, r}_N(M)$ is homotopic to $id_M$ rel $N$, then $h$ is $C^r$-isotopic to $id_M$ rel $N$ and 
$h \in {\cal D}^{\, r}_N(M)_0$.
\item[{\rm (2)}] Suppose $N$ is a compact smooth $2$-submanifold of $M$ with $\partial M \subset N$. 
\begin{itemize}
\item[{\rm (i)}] If $h \in {\cal D}^{\, r}_N(M)$ is $C^0$-isotopic to $id_M$ rel $N$, then $h$ is $C^r$-isotopic to $id_M$ rel $N$ and  $h \in {\cal D}^{\, r}_N(M)_0$. 
\item[{\rm (ii)}] $h \in {\cal D}^{\, r}_N(M)_0$ iff $h \in {\cal D}^{\, r}_N(M)$ and $h$ is $C^r$-isotopic to $id_M$ rel $N$ $($cf.\,Remark 2.1\,(2)$)$. 
\end{itemize}
\end{itemize}
\end{lemma}

In Corollary~\ref{cor_bdle} we have a principal bundle with fiber ${\mathcal G} \equiv {\cal D}^{\, r}_X(M)_0 \cap {\cal D}^{\, r}_N(M)$. 
The next theorem provides us with a sufficient condition which implies ${\mathcal G} = {\cal D}^{\, r}_N(M)_0$. 
The symbol $\#X$ denotes the number of elements (or cardinal) of a set $X$.

\begin{theorem}\label{thm_rel-iso}
Suppose $M$ is a compact connected smooth $2$-manifold $($possibly with boundary$)$, 
$N$ is a compact smooth $2$-submanifold of $M$ with $\partial M \subset N$ and $X$ is a subset of $N$.
Suppose $(M, N, X)$ satisfies the following conditions: 
\begin{itemize}
\item[(i)\ ] $M \neq {\Bbb T}$, ${\Bbb P}$, ${\Bbb K}$ or $X \neq \emptyset$. 
\item[(ii)\,] {\rm (a)} if $H$ is a disk component of $N$, then $\# (H \cap X) \geq 2$, \\
{\rm (b)} if $H$ is an annulus or M\"obius band component of $N$, then $H \cap X \neq \emptyset$, 
\item[(iii)]
{\rm (a)} if $L$ is a disk component of $cl(M \setminus N)$, then $\# (L \cap X) \geq 2$, \\ 
{\rm (b)} if $L$ is a M\"obius band component of $cl(M \setminus N)$, then $L \cap X \neq \emptyset$. 
\end{itemize}
Then we have 
\begin{enumerate}
\item if $h \in {\cal D}^{\, r}_N(M)$ is $C^0$-isotopic to $id_M$ rel $X$, then $h$ is $C^r$-isotopic to $id_M$ rel $N$, 
\item ${\cal D}^{\, r}_X(M)_0 \cap {\cal D}^{\, r}_N(M) = {\cal D}^{\, r}_N(M)_0$. 
\end{enumerate}
\end{theorem}

\begin{lemma}\label{lem_rel-iso} 
Suppose $M$ is a connected smooth $2$-manifold without boundary and $N$ is a smooth closed $2$-submanifold of $M$.
Suppose $N \neq \emptyset$, $cl(M \setminus N)$ is compact and $(M, N)$ satisfies the following conditions: 
\begin{itemize}
\item[(i)\ ] $M \not \cong {\Bbb T}$, ${\Bbb P}$, ${\Bbb K}$, 
\item[(ii)\,] each component $C$ of $\partial N$ does not bound a disk or a M\"obius band, 
\item[(iii)] each component of $N$ \, $\not\cong$ \, ${\Bbb S}^1 \times [0, 1]$, ${\Bbb S}^1 \times [0, 1)$.  
\end{itemize}
If $h \in {\cal D}^{\, r}_N(M)$ is $C^0$-isotopic to $id_M$, then 
$h$ is $C^r$-isotopic to $id_M$ rel $N$ and $h \in {\cal D}^{\, r}_N(M)_0$. 
\end{lemma} 

Theorem~\ref{thm_rel-iso} and Lemma~\ref{lem_rel-iso} follow from 
Lemma~\ref{lem_C^r-iso} and the corresponding statements in the $C^0$-case 
\cite[Theorem 3.1, Lemma 3.4 (and a remark after Lemma 3.4)]{Ya2}. 

\subsection{Diffeomorphism groups of noncompact $2$-manifolds} \mbox{} 

In this section we prove Theorems~\ref{thm_l2} and \ref{thm_HD}. 
Suppose $M$ is a noncompact connected smooth $2$-manifold without boundary and $X$ is a compact smooth submanifold of $M$ (i.e., 
a disjoint union of a compact smooth 2-submanifold of $M$ and finitely many smooth circles and points). 
We need to separate the next two cases:  
\begin{itemize}
\item[(I)\,] $(M, X) =$ (a plane, $\emptyset$), (a plane, 1\,pt), (an open M\"obius band, $\emptyset$) or (an open annulus, $\emptyset$). 
\item[(II)] $(M, X)$ is not Case (I). 
\end{itemize} 
Theorems~\ref{thm_l2} and \ref{thm_HD} are rewritten as follows: 

\begin{theorem}\label{thm_main_2} {\rm (1)} ${\cal D}^{\, r}_X(M)_0$ is an $\ell_2$-manifold. 
\begin{itemize}
\item[{\rm (2)}] ${\cal D}^{\, r}_X(M)_0 \simeq {\Bbb S}^1$ in Case {\rm (I)} \ and \ ${\cal D}^{\, r}_X(M)_0 \simeq \ast$ in Case {\rm (II)}. 
\item[(3)] ${{\cal D}_X^{\, r}(M)^c_0}^\ast$ is HD in ${\cal D}^{\, r}_X(M)_0$. 
\end{itemize}
\end{theorem}

\begin{proof} 
{Case (I)}: \ From the assumption it follows that $(M, X) = ({\rm Int}\,N, X)$, where $(N, X) =  ({\Bbb D}, \emptyset)$, $({\Bbb D}, 1pt)$, $({\Bbb M}, \emptyset)$ or $({\Bbb A}, \emptyset)$. 
By Theorem~\ref{thm_h-type_cpt} it is seen that ${\cal D}^{\, r}_X(N)_0 \simeq {\Bbb S}^1$ and 
\[ \mbox{${\cal D}^{\, r}_{L \times \{ 0 \}}(L \times [0,1]) = {\cal D}^{\, r}_{L \times \{ 0 \}}(L \times [0,1])_0 \simeq \ast$ \ \ and \ \ 
${\cal D}^{\, r}_{L \times \{ 0, 1 \}}(L \times [0,1])_0 \simeq \ast$} \]  
for any closed smooth 1-manifold $L$ (i.e., $L$ is a disjoint union of finitely many circles).  
Therefore, from Proposition~\ref{prop_criterion_cpt} for $n = 2$ follows  
the assertions (1), (3) and (2)\,(I) ${\cal D}^{\, r}_X(M)_0 \cong {\cal D}^{\, r}_X(N)_0 \simeq {\Bbb S}^1$. 
\vskip 1mm 
{Case (II)}: \ We can write as $M = \cup_{i=0}^{\infty} \, M_i$, where $M_0 = X$ and for each $i \geq 1$
\begin{itemize}
\item[(a)] $M_i$ is a nonempty compact connected smooth 2-submanifold of $M$ and $M_{i-1} \subset {\rm Int}\, M_i$, 
\item[(b)] for each connected component $L$ of $cl \, (M \setminus M_i)$, $L$ is noncompact and $L \cap M_{i+1}$ is connected. 
\end{itemize}
Note that $M$ is a plane (an open M\"obius band, an open annulus) 
iff infinitely many $M_i$'s are disks (M\"obius bands, annuli respectively), 
and that any subsequence of $M_i$ ($i \geq 1$) also satisfies the conditions (a) and (b). 
Thus, passing to a subsequence, we may assume that

\begin{itemize}
\item[(i)\ ] if $M$ is a plane, then each $M_i$ is a disk, 
\item[(ii)\,] if $M$ is an open M\"obius band, then each $M_i$ is a M\"obius band, 
\item[(iii)] if $M$ is an open annulus, then each $M_i$ is an annulus (and the inclusion $M_i \subset M_{i+1}$ is essential),  
\item[(iv)] if $M$ is not a plane, an open M\"obius band or an open annulus, then each $M_i$ is not a disk, an annulus or a M\"obius band.
\end{itemize} 

\noindent For each $i \geq 1$ let $U_i = {\rm int}\,M_i$.  

To apply Proposition~\ref{prop_criterion_diff} we have to verify the following conditions: 
\begin{itemize}
\item[$(\natural)_1$] ${\cal D}^{\, r}_{M_k\cup (M - U_i)}(M)_0 \simeq \ast$ \ for each $i > k \geq 0$, 
\item[$(\natural)_2$] ${\cal G}^i_{k,j} \equiv {\cal D}^{\, r}_{M_k\cup (M - U_j)}(M)_0 \cap {\cal D}^{\, r}_{M_i}(M) 
= {\cal D}^{\, r}_{M_i\cup (M - U_j)}(M)_0$ \ for each $j > i > k \geq 0$. 
\end{itemize}

Choose a small closed collar $E_i$ of $\partial M_i$ in $U_{i+1} \setminus U_i$ and set $M_i' = M_i \cup E_i \subset U_{i+1}$. 

$(\natural)_1$: ${\cal D}^{\, r}_{M_k\cup (M - U_i)}(M)_0 \cong 
{\cal D}^{\, r}_{M_k \cup E_i}(M_i') \simeq \ast$ by 
Theorem~\ref{thm_h-type_cpt}. 
\vskip 1mm 
$(\natural)_2$: 
It suffices to show that \ 
$${\cal D}^{\, r}_{M_k \cup E_j}(M_j')_0 \cap {\cal D}^{\, r}_{M_i \cup E_j}(M_j') = {\cal D}^{\, r}_{M_i \cup E_j}(M_j')_0.$$  
We apply Theorem~\ref{thm_rel-iso} to $(\tilde{M}, \tilde{N}, \tilde{X}) = (M_j', M_i \cup E_j, M_k \cup E_j)$.  
It remains to verify that this triple satisfies the conditions (i)--(iii) in Theorem~\ref{thm_rel-iso}: 
\vskip 1mm
(ii): The components of $\tilde{N}$ consists of $M_i$ and the annulus components of $E_j$. 
The latter components obviously satisfy the condition (ii)\,(b). 
By the choice of $\{ M_j \}_{j \geq 1}$, if $M_i \cong {\Bbb D}$ (${\Bbb M}$, ${\Bbb A}$), then 
$M$ is a plane (an open M\"obius band, an open annulus). 
Thus the assumption of Case (II) implies that 
$(M_i, X) \neq ({\Bbb D}, \emptyset), ({\Bbb D}, 1pt), ({\Bbb M}, \emptyset), ({\Bbb A}, \emptyset)$  
and that $M_i$ satisfies the condition (ii). 
\vskip 1mm
(iii): Each component $L$ of $cl(\tilde{M} \setminus \tilde{N}) = M_j \setminus U_i$ meets both $\partial M_i$ and $\partial M_j$. 
Thus $\partial L$ is not connected and $L \not\cong {\Bbb D}, {\Bbb M}$. 
This completes the proof. 
\end{proof}

%%%%%%%%%%%%%%%%%%%%%%%%%%%

\begin{proof}[\bf Proof of Corollary~\ref{cor_HD}] 
Since ${{\cal D}_X^{\, r}(M)^c_0}^\ast \subset {\cal D}_X^{\, r}(M)^c_0 \subset {\cal D}^{\, r}_X(M)_0$, 
the subgroup ${\cal D}_X^{\, r}(M)^c_0$ is also HD in ${\cal D}^{\, r}_X(M)_0$. Thus, the assertion (1) follows from Theorem~\ref{thm_l2}\,(1) and 
Lemma~\ref{lem_HD}\,(1)(ii), while 
the assertion (2) follows from Lemma~\ref{lem_HD}\,(1)(i).
\end{proof} 

\begin{proof}[\bf Proof of Proposition~\ref{prop_c*=c}] \mbox{} 
For notational simplicity we set $G = {\cal D}_X^{\, r}(M)$. 
We choose an exhausting sequence $\{ M_i \}_{i \geq 1}$ of $(M, X)$ as in the proof of Theorem~\ref{thm_main_2} Case (II). 

\noindent {\bf [A]} If $(M, X)$ satisfies the condition (a) or (b), then $G^c_0 = G^{c\,\ast}_0$: 

Given $h \in G^c_0$, there exists $i \geq 1$ such that $h \in G(U_i)$. 
We show that $h \in G(U_i)_0 \subset G^{c\,\ast}_0$. 
If $(M, X)$ satisfies the condition (b), then by Theorem~\ref{thm_h-type_cpt}\,(iv) 
$G(U_i) \cong {\cal D}^{\, r}_{E_i}(M_i \cup E_i) \simeq {\cal D}^{\, r}_{\partial M_i}(M_i) \simeq \ast$, where 
$E_i$ is a small closed collar of $\partial M_i$ in $M \setminus U_i$.  
This means that $G(U_i) = G(U_i)_0$. 
Below we assume that $(M, X)$ does not satisfy the condition (b). 
Note that $h$ is $C^0$-isotopic to $id_M$ rel $X$ since $h \in G_0$ and $G_0$ is a connected ANR. 
\vskip 1.5mm 

\noindent (1) The case (a) with $X = \emptyset$. 
\vskip 1mm 
Let $N = M \setminus U_i$ and we apply Lemma~\ref{lem_rel-iso} to $(M, N)$ and $h$. 
The conditions (ii) and (iii) in Lemma~\ref{lem_rel-iso} are verified as follows: 

(ii) If a circle component $C$ of $\partial N = \partial M_i$ bounds a disk or a M\"obius band $D$, then $M_i = D$ and by the choice of the exhausting sequence $\{ M_i \}_{i \geq 1}$ it follows that $M$ is a plane or an open M\"obius band. 
This contradicts the assumption that $(M, X)$ does not satisfy the condition (b)

(iii) Any component of $N$ is noncompact and not diffeomorphic to ${\Bbb S}^1 \times [0, 1)$ by the condition (a). 

Since $h \in G_0$ is $C^0$-isotopic to $id_M$, from Lemma~\ref{lem_rel-iso} it follows that $h$ is $C^r$-isotopic to $id_M$ rel $N$ and $h \in G(U_i)_0$. 
\vskip 1mm

\noindent (2) The case (a) with $X \neq \emptyset$ : 
\vskip 0.5mm 
Let $F_k$ ($k = 1, \cdots, m$) denote  
the connected components of the compact 2-submanifold $cl(M_i \setminus X)$ 
which are disks or M\"obius bands 
let $H_\ell$ ($\ell = 1, \cdots, n$) denote the remaining components. 
Set $N = X \cup (\cup_k F_k) \cup (M \setminus U_i)$. 
There exists a $C^0$-isotopy $h_t$ rel $X$ from $h$ to $id_M$. 
Since $\partial F_k \subset X$, it follows that 
$h_t$ maps each $F_k$ onto itself. 
Thus, if we define $h' \in G_N$ by $h' = h$ on $X \cup (\cup_\ell H_\ell) \cup (M \setminus U_i)$ and $h' = id$ on $\cup_k F_k$, 
then $h'$ is $C^0$-isotopic to $h$ rel $X \cup (M \setminus U_i)$.
and $C^0$-isotopic to $id_M$ rel $X$. 

Define $A \subset X$ by choosing two interior points from each component of $X$. 
Set $M' = M \setminus A$ and $N' = N \setminus A$. 
Then $h'|_{M'} \in {\cal D}^{\, r}_{N'}(M')$ is $C^0$-isotopic to $id_{M'}$. 
We show that $(M', N')$ satisfies the conditions (i)--(iii) in Lemma ~\ref{lem_rel-iso}. 
First note that $M'$ is noncompact and $L \equiv cl(M' \setminus N') = cl(M \setminus N) = \cup_\ell H_\ell$ is compact. 

(ii) Suppose a circle component $C$ of $\partial N = \partial M_i \cup (\partial X \setminus (\cup_k \partial F_k))$ bounds a disk or a M\"obius band $D$ in $M'$. 
Then $M$ is the union of $D$ and a connected submanifold $W = cl(M - D)$ with $\partial W = C$, and it follows that $N \subset W$. In fact, 
($\alpha$) $X \subset W$ since $A \subset W$ and each component of $X$ meets $A$, 
($\beta$) each $F_k \subset W$ since $F_k$ meets $X$ but does not meet $C$ and 
($\gamma$) $M \setminus U_i \subset W$ since each component of $M \setminus U_i$ is noncompact. 
If $C \subset \partial X \setminus (\cup_k \partial F_k)$, Thus $D \subset \cup_\ell H_\ell$ and $D = H_\ell$ for some $\ell$, but this contradicts the choice of $H_\ell$. If $C \subset \partial M_i$, then $D \supset M_i \supset A$, which  
contradicts that $D \subset M'$. 

(iii) Let $J$ be any component of $N'$. 
If $J \subset (X \setminus A) \cup (\cup_k F_k)$, 
then $J$ has two ends, since it is the union of $X_1 \setminus A$ for some component $X_1$ of $X$ and some $F_k$'s. 
Otherwise, $J$ is a component of $M \setminus U_i$, so it is noncompact and $J \not\cong {\Bbb S}^1 \times [0, 1)$ by the condition (a). 

Therefore, by Lemma~\ref{lem_rel-iso} $h'|_{M'}$ is $C^0$-isotopic to $id_{M'}$ rel $N'$ and by the end compactification 
we obtain a $C^0$-isotopy $h \simeq h' \simeq id_M$ rel $X \cup (M \setminus U_i)$. 
This implies that $h \in G(U_i)_0$.
\vskip 1mm

\noindent {\bf [B]} If $(M, X)$ does not satisfy the conditions (a) and (b), then $G^{c\,\ast}_0 \subsetneqq G^c_0$:
\vskip 1mm
Suppose $(M, X)$ does not satisfy the conditions (a) and (b). 
Then $M$ contains a product end $E \cong {\Bbb S}^1 \times [0, 1)$ such that $E \cap X = \emptyset$. 
We identify $E$ with ${\Bbb S}^1 \times [0, 1)$. 
Note that $F = cl(M \setminus E)$ is a connected submanifold of $M$ and $\partial F = \partial E$.
Let $h \in G^c$ denote a Dehn twist in $E$.  
We show that $h \in G^c_0 \setminus G^{c\,\ast}_0$. 
A $C^r$-isotopy $h_t : h \simeq id_M$ with $h_t \in G^c_0$ is obtained by 
sliding the Dehn twist towards $\infty$. This implies that $h \in G^c_0$.
It remains to show that $h \not\in G^{c\,\ast}_0$. 
On the contrary, suppose $h \in G^{c\,\ast}_0$. 
Then $h \in G(M_i)_0$ for some $i \geq 1$ and there exists a $C^0$-isotopy $h_t : h \simeq id_M$ rel $X$ with ${\rm supp}\,h_t \subset M_i$ ($0 \leq t \leq 1$).  

If $F \cong {\Bbb S}^1 \times [0, 1)$, then $M \cong {\Bbb S}^1 \times (-1, 1)$ and $M_i$ is contained in an annulus $L$ in $M$. 
This means that the Dehn twist in $L$ is $C^0$-isotopic to $id_L$ rel $\partial L$, a contradiction. 
Thus we may assume that $F \not\cong {\Bbb S}^1 \times [0, 1)$.

Choose $t \in (0, 1)$ such that $M_i \subset F \cup L$, where $L = {\Bbb S}^1 \times [0, t]$ is an annulus. 
Let $Y = {\Bbb S}^1 \times [t, 1)$ ($= cl(E \setminus L)$) and $N = F \cup Y$. 
Then $h|_N = id$ and $h_t : h \simeq id_M$ is a $C^0$-isotopy rel $X \cup Y$. 
Define $A$ by choosing an interior point from $Y$ and two interior points from $X$ if $X \neq \emptyset$. 
Let $M' = M \setminus A$ and $N' = N \setminus A$.
We show that $(M', N')$ satisfies the conditions (i)--(iii) in Lemma ~\ref{lem_rel-iso}. 

(ii) $\partial N' = \partial N$ consists of two circles $C_1 = {\Bbb S}^1 \times \{ 0 \}$ and $C_2 = {\Bbb S}^1 \times \{ t \}$. 
If $C_1$ bounds a disk or a M\"obius band $D$ in $M'$, then $F = D \subset M'$, 
so we have $X = \emptyset$ and $M = D \cup E$ is a plane or an open M\"obius band (i.e. the condition (b)). 
This contradicts the assumption in [B]. 
Similarly, the circle $C_2$ bound neither a disk nor a M\"obius band in $M'$.

(iii) Note that $N'$ consists of the two components $F \setminus A$ and $Y \setminus A$. 
These are homeomorphic to 
neither ${\Bbb S}^1 \times [0, 1]$ nor ${\Bbb S}^1 \times [0, 1)$.
In fact, 
$Y \setminus A$ has two ends, and if $X \neq \emptyset$, then 
$F \setminus A$ has at least two ends. 
If $X = \emptyset$, then $F \setminus A = F \not\cong  {\Bbb S}^1 \times [0, 1)$ and $\partial F$ is a circle. 

Since $h' = h|_{M'}$ is $C^0$-isotopic to $id_{M'}$ and $h'|_{N'} = id$, 
by Lemma~\ref{lem_rel-iso} we obtain a $C^0$-isotopy $h' \simeq id_{M'}$ rel $N'$. 
This means that 
the Dehn twist $h|L$ on the annulus $L$ is isotopic to $id_L$ rel $\partial L$. 
This is a contradiction. 
This completes the proof. 
\end{proof}

%%%%%%%%%%%%%%%%%%

\section{Groups of volume-preserving diffeomorphisms of noncompact 2-manifolds} 

In this final section we discuss topological types of groups of volume-preserving diffeomorphisms of noncompact 2-manifolds with the compact-open $C^\infty$-topology. Since we are only concerned with the groups of $C^\infty$-diffeomorphisms with the compact-open $C^\infty$-topology, we omit the symbol $\infty$ from the notations. 

\subsection{General properties of groups of volume-preserving diffeomorphisms of $n$-manifolds} \mbox{} \hspace*{2mm} 
Suppose $M$ is a connected oriented smooth $n$-manifold possibly with boundary,  
$X \subsetneqq M$ is a closed subset of $M$ and 
$\omega$ is a positive volume form on $M$. 
Let ${\cal D}_X(M; \omega) \subset {\cal D}^+_X(M)$ denote the subgroups of ${\cal D}_X(M)$ consisting of 
$\omega$-preserving diffeomorphisms and 
orientation-preserving diffeomorphisms of $M$ respectively, that is, 
$${\cal D}_X(M; \omega) = \big\{ h \in {\cal D}_X(M) \mid h^\ast \omega = \omega \big\}.$$ 
These subgroups are endowed with the subspace topology (i.e., the compact-open $C^\infty$-topology). 
As before, 
the subscript `0' denotes the identity connected component. 
Note that ${\cal D}^+_X(M)_0 = {\cal D}_X(M)_0$. 
The group ${\cal D}_X(M; \omega)$ is a separable, completely metrizable topological group 
since it is a closed subgroup of ${\cal D}_X(M)$. 
It is also seen to be infinite-dimensional and non locally compact. 
Hence, by Theorem~\ref{thm_top-gp} the group ${\cal D}_X(M; \omega)$ (or ${\cal D}_X(M; \omega)_0$) is an $\ell_2$-manifold iff it is an ANR.  

First we recall Moser's theorem \cite{Mos} and its extension to the noncompact  case \cite{Ya4}. (We refer to \cite{AP, Be1, Ya4} for end compactifications and related matters.) 
The space $E_M$ of ends of $M$ is a compact 0-dimensional metrizable space. 
Let $E_M^\omega$ denote the subspace of $E_M$ consisting of $\omega$-finite ends of $M$. Each $h \in {\cal H}(M)$ admits a unique homeomorphic extension $\overline{h}$ on the end compactification $M \cup E_M$. 
We define ${\cal D}^+(M, E_M^\omega) = \{ h \in {\cal D}^+(M) \mid \overline{h}(E_M^\omega) = E_M^\omega\}$. 

Let ${\cal V}^+(M; \omega(M), E_M^\omega)_{ew}$ denote the space of 
positive volume forms $\mu$ on $M$ such that $\mu(M) = \omega(M)$ and 
$E_M^\mu = E_M^\omega$. This space is endowed with the finite-ends weak $C^\infty$-topology $ew$ (cf. \cite{Be1, Ya4}).
The group ${\cal D}^+(M, E_M^\omega)$ acts on the space ${\cal V}^+(M; \omega(M), E_M^\omega)$ by the push-forward of forms and 
induces the orbit map at $\omega$, 
$$\pi_\omega : {\cal D}^+(M, E_M^\omega) \longrightarrow {\cal V}^+(M; \omega(M), E_M^\omega)_{ew}.$$ 

Moser's theorem \cite{Mos} and \cite[Corollary 1.1]{Ya4} assert that this orbit map admits a section into ${\cal D}_\partial(M)_0$. This implies the following relation between the groups ${\cal D}_X(M; \omega) \subset {\cal D}_X(M)$. 

\begin{proposition}\label{prop_diff_vp_SDR} \mbox{} 
\begin{itemize}
\item[(1)] 
$\big({\cal D}^{\infty\,+}(M;E_M^\omega), {\cal D}(M; \omega)\big) \cong 
\big({\cal V}^+(M; \omega(M), E_M^\omega)_{ew}, \{ \omega \}\big) \times {\cal D}(M; \omega)$, \\
$({\cal D}(M)_0, {\cal D}(M; \omega)_0) \cong 
({\cal V}^+(M; \omega(M), E_M^\omega)_{ew}, \{ \omega \}) \times {\cal D}(M; \omega)_0.$

\item[(2)] Suppose $N$ is a compact $n$-submanifold of ${\rm Int}\,M$. 
\begin{itemize}
\item[(i)\,] ${\cal D}_{M - N}(M; \omega)$ is a SDR of ${\cal D}_{M - N}^{\ +}(M)$ and is an ANR $($since ${\cal D}_{M - N}(M)$ is an ANR$)$.    
\item[(ii)] ${\cal D}_N(M; \omega)$ is a SDR of ${\cal D}_N^+(M, E_M^\omega)$ and ${\cal D}_N(M; \omega)_0$ is a SDR of ${\cal D}_N(M)_0$.
\end{itemize} 
\end{itemize} 
\end{proposition} 
\noindent In the statement (2) we apply Moser's theorem \cite{Mos} and \cite[Corollary 1.1]{Ya4} to each $L \in {\cal C}(M - {\rm Int}\,N)$. 

Next we recall the definition of the end charge homomorphism introduced by S.~R.~Alpern and V.\,S.\,Prasad \cite{AP}. 
Suppose $M$ is a noncompact connected orientable smooth $n$-manifold possibly with boundary and $\omega$ is a volume form on $M$. 
An end charge of $M$ is a finitely additive signed measure on the algebra of clopen subsets of $E_M$. 
Let ${\cal S}(M)$ denote the topological linear space of all end charges of $M$ with the weak topology and 
let ${\cal S}(M; \omega)$ denote the linear subspace of ${\cal S}(M)$ consisting of end charges $c$ of $M$ with $c(E_M) = 0$ and $c |_{E_M^\omega} = 0$. 

Let ${\cal D}_{E_M}(M; \omega) = \{ h \in  {\cal D}(M; \omega) \mid \overline{h}|_{E_M} = id_{E_M}\}$. 
We have ${\cal D}_{E_M}(M; \omega)_0 =  {\cal D}(M; \omega)_0$. 
For each $h \in {\cal D}_{E_M}(M; \omega)$ an end charge $c^\omega(h) \in {\cal S}(M; \omega)$ is defined by 
\[ c^\omega(h)(E_C) = \omega(C - h(C)) - \omega(h(C) - C), \]
where $C$ is any $n$-submanifold of $M$ such that ${\rm Fr}_M C$ is compact 
and $E_C \subset E_M$ is the set of ends of $C$. 
The end charge homomorphism 
$$c^\omega : {\cal D}_{E_M}(M; \omega) \to {\cal S}(M; \omega) : h \longmapsto c^\omega(h)$$ 
is a continuous group homomorphism. 
Let $c^\omega_0 : {\cal D}(M; \omega)_0 \to {\cal S}(M; \omega)$ denote the 
restriction of $c^\omega$ to ${\cal D}(M; \omega)_0$. 

The kernels ${\rm ker}\,c^\omega$ and ${\rm ker}\,c^\omega_0$ are separable, non locally compact, completely metrizable topological groups. 
Hence, by Theorem~\ref{thm_top-gp} the group ${\rm ker}\,c^\omega$ (or ${\rm ker}\,c^\omega_0$) is an $\ell_2$-manifold iff it is an ANR.  
In \cite[Corollary 1.2]{Ya4} we have shown that the homomorphism $c^\omega$ has a continuous (non-homomorphic) section into ${\cal D}_\partial(M; \omega)_0$. 
This clarifies the relation between the groups ${\rm ker}\,c^\omega_0 \subset {\cal D}(M; \omega)_0$. 

\begin{proposition}\label{prop_ker_SDR} \mbox{} 
\begin{itemize}
\item[(1)] $({\cal D}_{E_M}(M; \omega), {\rm ker}\,c^\omega) 
\cong ({\cal S}(M; \omega), 0) \times {\rm ker}\,c^\omega$, \ 
$({\cal D}(M; \omega)_0, {\rm ker}\,c^\omega_0) 
\cong ({\cal S}(M; \omega), 0) \times {\rm ker}\,c^\omega_0$. 
\item[(2)] 
\begin{itemize}
\item[{\rm (i)}\,] ${\rm ker}\,c^\omega$ is a SDR of ${\cal D}_{E_M}(M; \omega)$. 
\item[{\rm (ii)}] ${\rm ker}\,c^\omega_0 = ({\rm ker}\,c^\omega)_0$ and it is a SDR  of ${\cal D}(M; \omega)_0$. 
\end{itemize}
\end{itemize}
\end{proposition}

%%%%%%%%%%%%%%%%%%

\subsection{The bundle theorem} \mbox{}

Next we obtain the bundle theorems for groups of volume-preserving diffeomorphisms. 
Suppose $M$ is a connected oriented smooth $n$-manifold without boundary,  
$\omega$ is a positive volume form on $M$, 
$N$ is a compact smooth $n$-submanifold of $M$, 
$X$ is a closed subset of $N$ and 
$N_0$ is smooth $n$-submanifold of $M$ such that $N \subset U_0 \equiv {\rm Int}\,N_0$. 

For notational simplicity we set  
$(G, H, F) = \big({\cal D}(M), {\cal D}(M; \omega), {\rm ker}\,c^\omega\big)$. 
Let ${\cal V}^+(M)$ denote 
the space of positive volume forms on $M$ 
endowed with the compact-open $C^\infty$-topology. 
The group $G_0$ acts continuously on ${\cal V}^+(M)$ from the right by the pullback
$\mu \cdot h = h^\ast \mu$.  
For any subset ${\cal F}$ of ${\cal E}^\infty_X(N, M)$ 
the symbol ${\cal F}^{co}$ denotes the space ${\cal F}$ endowed with the subspace topology ($=$ the compact-open $C^\infty$-topology). 
When $i_N \in {\cal F}^{co}$, the symbol ${\cal F}^{co}_0$ denotes 
the connected component of the inclusion $i_N : N \subset M$ in ${\cal F}^{co}$. 
Hence, the space ${\cal E}^H_X(N, U_0)$ carries the quotient topology, while 
${\cal E}^H_X(N, U_0)^{co}$ carries the compact-open $C^\infty$-topology.
Let ${\cal C}(Y)$ denote the set of connected components of a space $Y$. 

The extension theorems for the transformation groups $H$ and $F$ are summarized as follows: 

\begin{theorem}\label{thm_extension_vp}
{\rm ($H$)} There exists a neighborhood ${\cal U}$ of $i_N$ in ${\cal E}^H_X(N, U_0)^{co}$ and 
\begin{itemize}
\item[] a map $\phi : {\cal U} \to H(U_0)_0 \cap H_X$ such that $\phi(f)|_N = f$ $(f \in {\cal U})$ and $\phi(i_N) = id_M$. 

\item[{\rm ($F$)}] Suppose $\partial N_0$ is compact and $U$ is an open neighborhood of $N$ in $U_0$ such that 
$U \cap L$ is connected for each $L \in {\cal C}(N_0 - {\rm Int}\,N)$. Then 
there exists a neighborhood ${\cal U}$ of $i_N$ in 
${\cal E}^F_X(N, U_0)^{co}$ and 
a map $\phi : {\cal U} \to F(U)_0 \cap F_X$ such that 
$\phi(f)|_N = f$ $(f \in {\cal U})$ and $\phi(i_N) = id_M$. 
\end{itemize}
\end{theorem} 

\begin{proof}
In each case we may assume that $X = \emptyset$, since $\phi(f)|_X = i_X$ 
if $\phi(f)|_N = f$ and $f|_X = i_X$.  

($H$) Choose any compact smooth $n$-submanifold $N_1$ of $U_0$ with 
$N \subset U_1 \equiv {\rm Int}\,N_1$.  
Consider the subspace of ${\cal V}^+(M)$ defined by 
$${\cal V}^\omega_N(U_1) = \{ \mu \in {\cal V}^+(M) \mid \mu 
= \omega \text{ on } N \cup (M - U_1) \text{ and } 
\mu(L) = \omega(L) \text{ for each } L \in {\cal C}(N_0 - {\rm Int}\,N) \}.$$ 
Applying Moser's theorem \cite{Mos} or \cite[Corollary 1.1]{Ya4} to each $L \in {\cal C}(N_0 - {\rm Int}\,N)$ we obtain a map 
$$\eta : {\cal V}^\omega_N(U_1) \to G_N(U_0)_0$$ 
such that $\eta(\mu)^\ast\mu = \omega$ 
$( \mu \in {\cal V}^\omega_N(U_1))$ and 
$\eta(\omega) = id_M$. 

By Theorem~\ref{thm_extension} there exists a neighborhood ${\cal U}_1$ of $i_N$ in ${\cal E}^\infty(N, U_1)$ 
and a map 
$$\psi : {\cal U}_1 \to G(U_1)_0$$ 
such that 
$\psi(f)|_N = f$ $(f \in {\cal U}_1)$ and $\psi(i_N) = id_M$. 
Then ${\cal U} = {\cal U}_1 \cap {\cal E}^H(N, U_0)^{co}$ is a neighborhood of $i_N$ in ${\cal E}^H(N, U_0)^{co}$. 

The map $\psi$ induces a map 
$$\chi : {\cal U} \to {\cal V}^\omega_N(U_1), \quad 
\chi(f) = \psi(f)^\ast\omega.$$  
We verify that $\chi(f) \in {\cal V}^\omega_N(U_1)$ $(f \in {\cal U})$. 
Since $f \in {\cal E}^H(N, U_0)$ there exists a $h \in H(U_0)$ with $h|_N = f$. 
Since $\psi(f)|_N = f$ and $\psi(f) = id$ on $M - U_1$, we have 
$\chi(f)|_N = f^\ast\omega = (h^\ast\omega)|_N = \omega|_N$ and 
$\chi(f)|_{M - U_1} = \omega|_{M - U_1}$. 
Since $h^{-1}\psi(f) \in G_N(U_0)$, for each $L \in {\cal C}(N_0 - {\rm Int}\,N)$ 
it follows that $h^{-1}\psi(f)(L) = L$ and $\psi(f)(L) = h(L)$, and that 
$\chi(f)(L) = \omega(\psi(f)(L)) = \omega(h(L)) = \omega(L)$.  

Finally, the required map $\phi : {\cal U} \to H(U_0)_0$ is defined by 
$$\phi(f) = \psi(f)\eta(\chi(f)) \ \ (f \in {\cal U}).$$ 
By Proposition~\ref{prop_diff_vp_SDR}\,(2)(ii) we have $\phi(f) \in H(U_0) \cap G(U_0)_0 = H(U_0)_0$. 
\vskip 1mm 
($F$) 
By the assumption we can find 
a compact smooth $n$-submanifold $N_1$ of $M$ such that $N \subset {\rm Int}\,N_1$, $N_1 \subset U$ and $N_1 \cap L$ is connected for each $L \in {\cal C}(N_0 - {\rm Int}\,N)$. Let $U_1 = {\rm Int}\,N_1$ and $N_1^* = N_1 - {\rm Int}\,N$. 

Consider the subspace of ${\cal V}^+(M)$ defined by 
$${\cal V}^\omega_N(U_1) = \{ \mu \in {\cal V}^+(M) \mid \mu = \omega \text{ on } N \cup (M - U_1) \text{ and } 
\mu(K) = \omega(K) \text{ for each } K \in {\cal C}(N_1^*) \}.$$ 
Applying Moser's theorem \cite{Mos} 
to each $K \in {\cal C}(N_1^*)$, we obtain 
a map 
$$\eta : {\cal V}^\omega_N(U_1) \to G_N(U_1)_0$$ such that $\eta(\mu)^\ast\mu = \omega$ $( \mu \in {\cal V})$ and 
$\eta(\omega) = id_M$. 

By Theorem~\ref{thm_extension} there exists a neighborhood ${\cal U}_1$ of $i_N$ in ${\cal E}^\infty(N, U_1)$ 
and a map $$\psi : {\cal U}_1 \to G(U_1)_0$$ such that 
$\psi(f)|_N = f$ $(f \in {\cal U}_1)$ and $\psi(i_N) = id_M$. 
Then ${\cal U} = {\cal U}_1 \cap {\cal E}^F(N, U_0)^{co}$ is a neighborhood of 
$i_N$ in ${\cal E}^F(N, M)^{co}$.  

The map $\psi$ induces a map 
$$\chi : {\cal U} \to {\cal V}^\omega_N(U_1), \ \ \chi(f) = \psi(f)^\ast\omega.$$ 
We show that $\chi(f) \in {\cal V}^\omega_N(U_1)$ $(f \in {\cal U})$. 
Since $f \in {\cal E}^F(N, U_0)$ there exists an $h \in F(U_0)$ with $h|_N = f$. 
Since $\psi(f)|_N = f$ and $\psi(f) = id$ on $M - U_1$, we have 
$\chi(f)|_N = f^\ast\omega = (h^\ast\omega)|_N = \omega|_N$ and 
$\chi(f)|_{M - U_1} = \omega|_{M - U_1}$. 
By the choice of $N_1$, for each $K \in {\cal C}(N_1^*)$ there exists a unique $L \in {\cal C}(N_0 - {\rm Int}\,N)$ with $K = N_1 \cap L$. 
Since $h^{-1}\psi(f) \in G_N(U_0)$, it follows that $h^{-1}\psi(f)(L) = L$ and $\psi(f)(L) = h(L)$.  
Since $L$ is an $n$-submanifold of $M$ and ${\rm Fr}\,L = \partial L \subset \partial N \cup \partial N_0$ is compact, from the definition of the end charge $c^\omega_0(h)$ it follows that 
$$c^\omega_0(h)(E_L) = \omega(L \setminus h(L)) - \omega(h(L) \setminus L). $$ 
Since $c^\omega_0(h) = 0$, we have $\omega(L \setminus h(L)) = \omega(h(L) \setminus L)$. 
Let $L_1 = L - N_1$. Then, $L_1 = \psi(f)(L_1) \subset h(L)$ and it is seen that 
\begin{align*}
& K = L - L_1 = \big[ L \setminus h(L)\big] \cup \big[ (h(L) \cap L) \setminus L_1\big] \quad \text{and} \\
& \psi(f)(K) = \psi(f)(L - L_1) = \psi(f)(L) - L_1 = 
h(L) - L_1 = \big[ h(L) \setminus L\big] \cup \big[ (h(L) \cap L) \setminus L_1\big].
\end{align*}
Thus we have $\chi(f)(K) = \omega(\psi(f)(K)) = \omega(K)$ and this means that $\chi(f) \in {\cal V}^\omega_N(U_1)$. 

Since $H(U_1)$ is a SDR of $G(U_1)$, we have  
$H(U_1) \cap G(U_1)_0 = H(U_1)_0 \subset F(U)_0$. 
Therefore, the required map $\phi : {\cal U} \to F(U)_0$ is defined by 
\[ \mbox{$\phi(f) = \psi(f)\eta(\chi(f))$ \ $(f \in {\cal U}).$} \] 
\vskip -7mm
\end{proof}

\begin{remark} 
The statement ($F$) does not necessarily hold when $\partial N_0$ is non compact. 
An example is easily obtained by inspecting 
a case where two ends of $N_0$ is included in an end of $M$. 
\end{remark} 

Consider the restriction map $\pi : F_X(U) \to {\cal E}^F_X(N, U)^{co}$, $\pi(h) = h|_N$. 
The group $F_N(U)$ acts on $F_X(U)$ by the right translation. 
The next corollary easily follows form Theorem~\ref{thm_extension_vp}\,(F) and  Lemma~\ref{lem_loc-sec}. 

\begin{corollary}\label{cor_bdle_vp} 
Suppose $\partial N_0$ is compact and 
$U$ is an open neighborhood of $N$ in $U_0$ such that 
$U \cap L$ is connected for each $L \in {\cal C}(N_0 - {\rm Int}\,N)$. 
Then, the following hold. 
\begin{itemize}
\item[{\rm (1)}] The subspace ${\cal E}^F_X(N, U)^{co}$ is open in ${\cal E}^F_X(N, U_0)^{co}$ and the restriction map $$\pi : F_X(U) \to {\cal E}^F_X(N, U)^{co}$$ is a principal bundle with fiber $F_N(U)$. 

\item[{\rm (2)}] If $F_X(U)_0$ is open in $F_X(U)$, 
then the subspace ${\cal E}^F_X(N, U)^{co}_0$ is 
closed and open in ${\cal E}^F_X(N, U)^{co}$ and the  restriction map 
$$\pi : F_X(U)_0 \to {\cal E}^F_X(N, U)^{co}_0$$ 
is a principal bundle with fiber $F_X(U)_0 \cap F_N$. 
\end{itemize}
\end{corollary}

We also need the following complementary resutls to Theorem~\ref{thm_extension_vp}. 

\begin{lemma}\label{lem_ext_vp} Suppose $\partial N_0$ is compact and 
$F_X(N_1)_0$ is open in $F_X(N_1)$ for any compact $n$-submanifold $N_1$ of $U_0$ with $X \subset {\rm Int}\,N_1$.  

\begin{itemize} 
\item[(1)] Suppose $f \in {\cal E}^F_X(N, U_0)^{co}$ and 
$U$ is an open neighborhood of $f(N)$ in $U_0$ such that 
$U \cap L$ is connected for each $L \in {\cal C}(N_0 - {\rm Int}\,f(N))$. 
Then there exists a neighborhood ${\cal V}$ of $f$ in 
${\cal E}^F_X(N, U_0)^{co}$ and 
a map $\psi : {\cal V} \to F_X(U)_0$ such that 
$\psi(g)f= g$ $(g \in {\cal V})$ and $\psi(f) = id_M$. 

\item[(2)] Each $f \in {\cal E}^F_X(N, U_0)_0^{co}$ satisfies the following condition. 
\begin{itemize} 
\item[$(\ast)$] There exists 
a compact $n$-submanifold $N_1$ of $U_0$ such that 
${\cal E}^F_X(N, N_1)^{co}_0$ is a neighborhood of $f$ in ${\cal E}^F_X(N, U_0)_0^{co}$. 
\end{itemize}
\end{itemize} 
\end{lemma} 

\begin{proof} 
(1) We may assume that $U = {\rm Int}\,N_1$ for some 
compact $n$-submanifold $N_1$ of $U_0$.  
Applying Theorem~\ref{thm_extension_vp}\,(F) to $(M, N_0, U, f(N), X)$, we obtain a neighborhood ${\cal U}$ of $i_{f(N)}$ in ${\cal E}^F_X(f(N), U_0)^{co}$ and  
a map $\phi : {\cal U} \to F(U)_0 \cap F_X$ such that 
$\phi(k)|_{f(N)} = k$ $(k \in {\cal U})$ and $\phi(i_{f(N)}) = id_M$. 
By the assumption $F_X(U)_0$ is open in $F_X(U)$ and so 
we may assume that $\phi({\cal U}) \subset F_X(U)_0$. 

Consider the homeomorphism 
\[ \mbox{$\chi : {\cal E}^F_X(f(N), U_0)^{co} \approx {\cal E}^F_X(N, U_0)^{co}$, \ \ $\chi(k) = kf$.} \]  
Then ${\cal V} = \chi({\cal U})$ and 
$\psi = \phi\chi^{-1} : {\cal V} \to F_X(U)_0$ 
satisfy the required conditions. 

(2) We have to show that the subset 
${\cal F} = \big\{ f \in {\cal E}^F_X(N, U_0)^{co}_0 \ \big| \ (\ast) \big\}$ 
 coincides with ${\cal E}^F_X(N, U_0)^{co}_0$.  
Corollary~\ref{cor_bdle_vp} implies that $i_N \in {\cal F}$. 
Therefore, it suffices to show that 
${\cal F}$ is closed and open in ${\cal E}^F_X(N, N_1)^{co}_0$.  
From the definition itself we see that ${\cal F}$ is open in ${\cal E}^F_X(N, N_1)^{co}_0$. 

To see that ${\cal F}$ is closed, take any $f \in cl\,{\cal F}$. 
There exists a compact $n$-submanifold $N_f$ of $U_0$ such that 
$f(N) \subset U_f \equiv {\rm Int}\,N_f$ and 
$U_f \cap L$ is connected for each $L \in {\cal C}(N_0 - {\rm Int}\,f(N))$. 
By (1) we obtain a neighborhood ${\cal V}$ of $f$ in 
${\cal E}^F_X(N, U_0)^{co}_0$ and 
a map $\psi : {\cal V} \to F_X(U_f)_0$ such that 
$\psi(g)f = g$ $(g \in {\cal V})$ and $\psi(f) = id_M$. 
Since ${\cal V} \cap {\cal F} \neq \emptyset$, 
we can choose a $g \in {\cal V} \cap {\cal F}$. 
Since $g \in {\cal V}$ it follows that $\psi(g)f = g$ and $f = \psi(g)^{-1}g$.  
In turn, since $g \in {\cal F}$, there exists 
a compact $n$-submanifold $N_g$ of $U_0$ such that $N \subset U_g \equiv {\rm Int}\,N_g$ and  $g \in {\cal E}^F_X(N, N_g)^{co}_0$. 
We may assume that $U_g \cap L$ is connected for each $L \in {\cal C}(N_0 - {\rm Int}\,N)$. 
Then, by Corollary~\ref{cor_bdle_vp} 
the  restriction map 
$$\pi : F_X(U_g)_0 \to {\cal E}^F_X(N, U_g)^{co}_0$$  
is a principal bundle.  
Hence, there exists an $h \in F_X(U_g)_0$ such that $g = h|_N$.  

Take a compact $n$-submanifold $N_1$ of $U_0$ such that 
$N_f \cup N_g \subset N_1$. 
Then, we have ${\cal V} \subset {\cal E}^F_X(N, N_1)^{co}_0$.  
In fact, for any $k \in {\cal V}$, 
it follows that $\psi(k)\psi(g)^{-1}h \in (F_X(U_f)_0)^2 F_X(U_g)_0 \subset F_X(N_1)_0$ and that 
$k = \psi(k)f = \psi(k)\psi(g)^{-1}h|_N \in {\cal E}^F_X(N, N_1)^{co}_0$. 
This means that $f \in {\cal F}$. 
This completes the proof.   
\end{proof} 

%%%%%%%%%%%%%%%%

\subsection{Groups of volume-preserving diffeomorphisms of noncompact $n$-manifolds} \mbox{} 

Suppose $M$ is a noncompact connected orientable smooth $n$-manifold without boundary, $\omega$ is a volume form on $M$ 
and $X$ is a compact smooth $n$-submanifold of $M$. 
We set $(G, H, F) = \big({\cal D}(M), {\cal D}(M; \omega), {\rm ker}\,c^\omega \big)$. 
Choose a smooth exhausting sequence $\{ M_i \}_{i \geq 0}$ of $M$ 
such that $M_0 = X$ and for each $i \geq 1$
\begin{itemize} 
\item[(a)] $M_i$ is connected,  
\item[(b)] $L$ is noncompact and $L \cap M_{i+1}$ is connected  
for each $L \in {\cal C}(M - {\rm Int}\,M_i)$.  
\end{itemize}
Let $U_i = {\rm Int}\,M_i$ $(i \geq 1)$. 

\begin{lemma}\label{lem_exh-seq}
The tuple $(F_X, M, \{ M_i \}_{i \geq 1})$ satisfies the assumption $(A)$. 
\end{lemma}

\begin{proof} 
(A-0) Since $(G, M)$ has a weak topology, so does $(F, M)$. 

(A-1) Corollary~\ref{cor_bdle_vp} implies the following conclusions. 
For each $j > i > k \geq 0$ 
\begin{itemize} 
\item[(a)] ${\cal E}^F_{M_k}(M_i, U_j)^{co} = {\cal E}^F_{M_k}(M_i, U_j)$ and 
${\cal E}^F_{M_k}(M_i, M)^{co} = {\cal E}^F_{M_k}(M_i, M)$, since 
the restriction maps 
$\pi^i_{k,j} : F_{M_k}(U_j) \longrightarrow {\cal E}^F_{M_k}(M_i, U_j)^{co}$ and $\pi^i_{k} : F_{M_k} \to {\cal E}^F_{M_k}(M_i, M)^{co}$ are principal bundles,  

\item[(b)] the restriction map 
$\pi^i_{k,j} : F_{M_k}(U_j)_0 \longrightarrow {\cal E}^F_{M_k}(M_i, U_j)_0$ 
is a principal bundle with the structure group $F_{M_k}(U_j)_0 \cap F_{M_i}$. 
\end{itemize}

(A-2)\,(i) $F_X(U_i) = H_X(U_i)$ is an ANR for each $i \geq 1$. 

(A-2)(ii) By (A-1)(a) we can work under the compact-open $C^\infty$-topology. 
Let ${\cal U}_{k,j}^i = {\cal E}^F_{M_k}(M_i, U_j)^{co}_0$ 
\begin{itemize} 
\item[] $(j > i > k \geq 0)$. 
\item[(a)] ${\cal U}_{k,j}^i$ is an open subspace of ${\cal E}^F_{M_k}(M_i, M)^{co}_0$ : 
This follows from Corollary~\ref{cor_bdle_vp}. 

\item[(b)] ${\cal E}^F_{M_k}(M_i, M)^{co}_0 = \cup_{j > i}\,{\cal U}_{k,j}^i$ : 
By Lemma~\ref{lem_ext_vp}\,(2)  
for each $f \in {\cal E}^F_{M_k}(M_i, M)^{co}_0$ 
there exists a compact $n$-submanifold $N_1$ of $M$ such that 
${\cal E}^F_{M_k}(M_i, N_1)^{co}_0$ is a neighborhood of $f$ in 
${\cal E}^F_{M_k}(M_i, M)_0^{co}$. 
If we choose a $j > i$ such that $N_1 \subset M_j$, then we have 
$f \in {\cal U}_{k,j}^i$. 

\item[(c)] $cl\,{\cal U}_{k,j}^i \subset {\cal U}_{k,j+1}^i$ : 
Given any $f \in cl\,{\cal U}_{k,j}^i$. 
Then $f(M_i) \subset M_j \subset U_{j+1}$.  
First we show that $L \cap U_{j+1}$ is connected 
for each $L \in {\cal C}(M - {\rm Int}\,f(M_i))$. 
Take $g \in {\cal U}_{k,j}^i$ which is sufficiently close to $f$ so that 
there exists a $k \in G(U_{j+1})$ such that $kg  = f$. 
Since the restriction map 
$\pi^i_{k,j} : F_{M_k}(U_j)_0 \longrightarrow {\cal U}_{k,j}^i$ is surjective, 
there exists an $h \in F_{M_k}(U_j)_0 \subset G(U_{j+1})$ such that 
$h|_{M_i} = g$. 
Recall the condition (b) for the exhausting sequence $\{ M_i \}_{i \geq 0}$. 
Then the claim is verified by the homeomorphism of tuples $kh : (M, U_{j+1}, M_i) \approx (M, U_{j+1}, f(M_i))$.  

By Lemma~\ref{lem_ext_vp}\,(1) there exists a neighborhood ${\cal V}$ of $f$ in 
${\cal E}^F_X(N, M)^{co}_0$ and 
a map $\psi : {\cal V} \to F_{M_k}(U_{j+1})_0$ such that 
$\psi(g)f= g$ $(g \in {\cal V})$ and $\psi(f) = id_M$. 
Since ${\cal V} \cap {\cal U}_{k,j}^i \neq \emptyset$, we can choose a $g' \in {\cal V} \cap {\cal U}_{k,j}^i$. 
There exists an $h' \in F_{M_k}(U_j)_0 \subset F_{M_k}(U_{j+1})_0$ such that 
$h'|_{M_i} = g'$. 
It follows that 
$\psi(g')^{-1}h' \in F_{M_k}(U_{j+1})_0$ and 
$f = \psi(g')^{-1}h'|_{M_i} \in {\cal U}_{k,j+1}^i$ as required. 
\end{itemize}
\vskip -7mm
\end{proof} 

\begin{remark} 
The tuple $(H_X, M, \{ M_i \}_{i \geq 1})$ does not satisfy the assumption (A-2)\,(ii). 
\end{remark}

Lemma~\ref{lem_(C-2)} and Theorem~\ref{thm_criterion}\,(2)(ii) induce the following conclusions. 

\begin{proposition}\label{prop_criterion_vp_hd}
{\rm (1)} If $F_X(U_j)_0 \cap F_{M_i} \simeq \ast$ for each $j > i \geq 1$, then 
$(F_X)^{c\,\ast}_0$ is HD in $(F_X)_0$. 
\begin{itemize} 
\item[(2)] If $F_{M_1}$ is connected and $F_{M_1}(U_j)_0 \cap F_{M_i} \simeq \ast$ for each $j > i \geq 2$, then 
$(F_X)^{c\,\ast}_0$ is HD in $(F_X)_0$. 
\end{itemize} 
\end{proposition} 

For $n$-manifolds of finite type, Proposition~\ref{prop_criterion_cpt} and 
Proposition~\ref{prop_criterion_vp_hd} 
induce the following conclusions. 

\begin{proposition}\label{prop_diff_vp_finite} 
Suppose $M = {\rm Int}\,N$ for some compact connected orientable $n$-manifold $N$ with non-empty boundary, $\omega$ is a volume form on $M$
and $X$ is a compact smooth $n$-submanifold of $M$. 
Then the following hold. 
\begin{itemize} 
\item[(1)] Both ${\cal D}_X(M; \omega)_0$ and 
$(F_X)_0$ are $\ell^2$-manifolds. 
\item[(2)] If ${\cal D}_{\partial N \times \{ 0, 1 \}}(\partial N \times [0,1])_0 \simeq \ast$, then $(F_X)^{c\,\ast}_0$ is HD in $(F_X)_0$. 
\end{itemize}
\end{proposition} 

\begin{proof} 
(1) The group $(G_X)_0$ is an ANR by Proposition~\ref{prop_criterion_cpt}.   
Hence, so is $(H_X)_0$ by Proposition~\ref{prop_diff_vp_SDR}.
We can apply Proposition~\ref{prop_diff_vp_SDR}\,(2)(ii) to each $L \in {\cal C}(M - {\rm Int}\,X)$ to show that $(F_X)_0$ is a SDR of $(H_X)_0$. 
Hence $(F_X)_0$ is also an ANR. 

(2) We construct an exhausting sequence $\{ M_i \}_{i \geq 1}$ of $(M, X)$ 
as in the proof of Proposition~\ref{prop_criterion_cpt}. 
Let $U_i = {\rm Int}\,M_i$ $(i \geq 1)$. 
Since this exhausting sequence satisfies the above conditions (a), (b),
the tuple $(F_X, M, \{ M_i \}_{i \geq 1})$ satisfies the assumption (A). 
Note that each $L \in {\cal C}(M - {\rm Int}\,M_1)$ is a product end. 

By Proposition~\ref{prop_diff_vp_finite}\,(2) it suffices to show the following statements. 

\begin{itemize}
\item[(i)\,] $F_{M_1}$ is connected :  
By Proposition~\ref{prop_diff_vp_SDR} 
the subgroup $H_{M_1}$ is a SDR of 
$G_{M_1} = {\cal D}_{M_1}^+(M, E_M^\omega)$.  
We can apply Proposition~\ref{prop_ker_SDR}\,(2)(i) to each $L \in {\cal C}(M - {\rm Int}\,M_1)$ to show that 
$F_{M_1}$ is a SDR of $H_{M_1}$. 
Since $G_{M_1} \simeq \ast$ by the Alexander trick towards $\infty$, 
we have $F_{M_1} \simeq \ast$. 

\item[(ii)] $F_{M_1}(U_j)_0 \cap F_{M_i} \simeq \ast$ for $j > i \geq 2$: 
In the proof of Proposition~\ref{prop_criterion_cpt}\,(3)  
we have already shown that 
$G_{M_1}(U_j)_0 \cap G_{M_i} \simeq \ast$ for $j > i \geq 2$.  
By Proposition~\ref{prop_diff_vp_SDR}\,(1) 
$F_{M_i}(U_j) = H_{M_i}(U_j)$ is a SDR of $G_{M_i}(U_j)$. 
Then, it follows that $F_{M_1}(U_j)_0 \cap F_{M_i} = F_{M_i}(U_j)_0$ and 
it is a SDR of 
$G_{M_1}(U_j)_0 \cap G_{M_i} = G_{M_i}(U_j)_0$. 
This implies the assertion. 
\end{itemize}
\vskip -7mm 
\end{proof} 

\subsection{Groups of volume-preserving diffeomorphisms of noncompact 2-manifolds} \mbox{} 

Suppose $M$ is a noncompact connected orientable smooth 2-manifold without boundary and $\omega$ is a volume form on $M$. 

\begin{proof}[\bf Proof of Theorem~\ref{thm_vol-pre}]  
Since ${\cal D}(M; \omega)_0$ is a SDR of ${\cal D}(M)_0$ (\cite[Corollary 1.1\,(2)(ii)]{Ya4}, see Section 6.1),  
the assertions follow from Theorem~\ref{thm_l2} and the observations in Section 6.1. 
\end{proof}

\begin{proof}[\bf Proof of Theorem~\ref{thm_vol-pre_hd}]  
(1) Since ${\rm ker}\,c^\omega_0$ is a SDR of ${\cal D}(M; \omega)_0$ (\cite[Corollary 1.2\,(2)(ii)]{Ya4}, see Section 6.1), 
the assertion follows from Theorem~\ref{thm_vol-pre} and the observations in Section 6.1. 

(2) Let $(G, H, F) = \big({\cal D}(M), {\cal D}(M; \omega), {\rm ker}\,c^\omega \big)$. 
We have to show that $F^{c\,\ast}_0$ is HD in $F_0$. 
We separate the next two cases (cf. Theorem~\ref{thm_main_2}):  
\begin{itemize}
\item[(I)\,] $M =$ a plane or an open annulus. \quad 
(II) $M$ is not Case (I). 
\end{itemize} 

{Case (I)}: \ Since $M = {\rm Int}\,N$ ($N =  {\Bbb D}$ or ${\Bbb A}$), the 
assertion follows from Proposition~\ref{prop_diff_vp_finite}\,(2) (cf. )  

{Case (II)}: \ There exists an exhausting sequence $\{M_i\}_{i \geq 1}$ of $M$ which satisfies the conditions (a), (b) in Section 6.3 and such that 
each $M_i$ is neither a disk, an annulus nor a M\"obius band. 
Let $U_i = {\rm int}\,M_i$ $(i \geq 1)$.  
In the proof of Theorem~\ref{thm_main_2}, we have shown that 
$$G(U_j)_0 \cap G_{M_i} = G_{M_i}(U_j)_0 \simeq \ast \ \ \text{ for each \ $j > i \geq 1$.}$$ 
Since $H_{M_i}(U_j)$ is a SDR of $G_{M_i}(U_j)$, this implies that  
$$F(U_j)_0 \cap F_{M_i} = H(U_j)_0 \cap H_{M_i} = H_{M_i}(U_j)_0 \simeq \ast \ \ \text{ for each \ $j > i \geq 1$.}$$ 
Hence, the conclusion follows from Proposition~\ref{prop_criterion_vp_hd}.
\end{proof}

%%%%%%%%%%%%%%%%%%


\begin{thebibliography}{99}

\bibitem{AP}
S.~R.~Alpern and V.\,S.\,Prasad, {Typical dynamics of volume-preserving homeomorphisms,} 
Cambridge Tracts in Mathematics, Cambridge University Press@(2001). 

\bibitem{BMSY} 
T. Banakh, K. Mine, K. Sakai and T. Yagasaki, 
 {\em Homeomorphism and diffeomorphism groups of non-compact manifolds with the Whitney topology}, preprint (arXiv math.GT/0802.0337). 

\bibitem{Bany1}
A.~Banyaga, {Formes-volume sur les vari\'et\'es \`a bord,} 
{\sl Enseignement Math.} (2) 20 (1974) 127 - 131. 

\bibitem{Bany2}
A.~Banyaga, {The Structure of Classical Diffeomorphism groups,} 
{\sl Mathematics and Its Applications} 400, Kluwer Academic Publishers Group, Dordrecht, 1997. 

\bibitem{Be1}
R.~Berlanga, {Groups of measure-preserving homeomorphisms as deformation retracts,} 
{\sl J. London Math. Soc.\,(2)} 68 (2003) 241 - 254. 

\bibitem{BP}
C. Bessaga and A. Pe{\l}czy\'{n}ski, 
Selected topics in infinite-dimensional topology, 
Polska Akademia Nauk Instytut Mate., Monografie Mate., 58, 
PWN Polish Scientific Publishers, Warszawa, 1975.

\bibitem{Ce}
J. Cerf, 
{\em Topologie de certains espaces de plongements}, 
Bull.\,Soc.\,Math.\,France, 89 (1961) 227--380. 

\bibitem{DT}
T. Dobrowolski and H. Toru\'nczyk, 
{\em Separable complete ANR's admitting a group structure are Hilbert manifolds}, Topology Appl., 12 (1981) 229--235. 

\bibitem{EE}
C. J. Earle and J. Eells, 
{\em A fiber bundle discription of Teichm\"uller theory}, 
J.~Diff. Geom., 3 (1969) 19--43.

\bibitem{Ep}
D. B. A. Epstein, 
{\em Curves on 2-manifolds and isotopies}, 
Acta Math., 155 (1966) 83--107.

\bibitem{Ge}
R. Geoghegan, 
{\em On spaces of homeomorphisms, embeddings, and functions - I}, 
Topology, 11 (1972) 159--177. 

\bibitem{GS}
R.~E.~Greene and K.~Shiohama, Diffeomorphisms and volume-preserving embeddings of noncompact manifolds, 
{\sl Trans.\,Amer.\,Math.\,Soc.,} 255 (1979) 403 - 414. 

\bibitem{Ham}
R. S. Hamilton, 
{\em The inverse function theorem of Nash and Moser}, 
Bull. Amer. Math. Soc. (New Series), 7 n.1 (1982) 65--222.

\bibitem{Ha}
M. E. Hamstrom, 
{\em Homotopy groups of the space of homeomorphisms on a 2-manifold},
Illinois J. Math., 10 (1966) 563 - 573. 

\bibitem{Han}
O. Hanner, 
{\em Some theorems on absolute neighborhood retracts}, Ark. Mat., 1 (1951) 389--408.

\bibitem{Hi}
M. W. Hirsch, 
Differential Topology, 
GTM 33, Springer-Verlag, New York, 1976.

\bibitem{Hu}
S. T. Hu, 
Theory of Retracts, 
Wayne State Univ.\,Press, Detroit, 1965.

\bibitem{KM}
J. Kalliongis and D. McCullough, 
{\em Fiber-preserving imbeddings and diffeomorphisms}, preprint. 

\bibitem{Les} J. A. Leslie, 
{\em On a differential structure for the group of diffeomorphisms}, 
Topology 6 (1967), 263--271.

\bibitem{LM}
R. Luke and W. K. Mason, 
{\em The space of homeomorphisms on a compact two - manifold is an absolute neighborhood retract},
Trans. Amer. Math. Soc., 164 (1972), 275 - 285. 

\bibitem{Mos} 
J.~Moser, On the volume elements on a manifold, 
{\sl Trans.\,Amer.\,Math.\,Soc.,} 120 (1965) 286 - 294. 

\bibitem{Pa1}
R. S. Palais, 
{\em Local triviality of the restriction map for embeddings}, 
Comment Math.~Helv., 34 (1960) 305--312.

\bibitem{Pa2}
\bysame, 
{\em Homotopy theory of infinite dimensional manifolds}, 
Topology, 5 (1966) 1--16.

\bibitem{Sc}
G. P. Scott, 
{\em The space of homeomorphisms of 2-manifold}, 
Topology, 9 (1970) 97--109.

\bibitem{Seeley}
R. T. Seeley, 
{\em Extension of $C\sp{\infty }$ functions defined in a half space},
Proc.\ Amer.\ Math.\ Soc.\ {\bf 15} (1964), 625--626.

\bibitem{Sm}
S. Smale, 
{\em Diffeomorphisms of the 2-sphere}, 
Proc. Amer. Math. Soc., 10 (1959) 621--626.

\bibitem{To}
H. Toru\'nczyk, 
{\em Characterizing Hilbert space topology}, 
Fund. Math., 111 (1981) 247--262.

\bibitem{Ya1}
T. Yagasaki,  
{\em Spaces of embeddings of compact polyhedra into 2-manifolds}, 
Topology Appl., 108 (2000) 107--122.

\bibitem{Ya2}
\bysame, 
{\em Homotopy types of homeomorphism groups of noncompact 2-manifolds}, Topology Appl., 108 (2000) 123--136.

\bibitem{Ya3}
\bysame, 
{\em The groups of PL and Lipschitz homeomorphisms of noncompact 2-manifolds}, 
Bulletin of the Polish Academy of Sciences, Mathematics, 51 (2003) 445--466.

\bibitem{Ya4}
\bysame, 
{\em Groups of volume-preserving diffeomorphisms of noncompact manifolds and mass flow toward ends}, 
preprint (arXiv math.GT/0805.3552). 
\end{thebibliography}
\end{document}